\def\Date{2008/05/09}


\ifx\pdfoutput\jamaisdefined\else
\input supp-pdf.tex \pdfoutput=1 \pdfcompresslevel=9

\fi

%

\magnification=1200
\hsize=11.25cm
\vsize=18cm
\parskip 0pt
\parindent=12pt
\voffset=1cm
\hoffset=1cm



\catcode'32=9

\font\tenpc=cmcsc10
\font\eightpc=cmcsc8
\font\eightrm=cmr8
\font\eighti=cmmi8
\font\eightsy=cmsy8
\font\eightbf=cmbx8
\font\eighttt=cmtt8
\font\eightit=cmti8
\font\eightsl=cmsl8
\font\sixrm=cmr6
\font\sixi=cmmi6
\font\sixsy=cmsy6
\font\sixbf=cmbx6

\skewchar\eighti='177 \skewchar\sixi='177
\skewchar\eightsy='60 \skewchar\sixsy='60

\catcode`@=11

\def\tenpoint{%
  \textfont0=\tenrm \scriptfont0=\sevenrm \scriptscriptfont0=\fiverm
  \def\rm{\fam\z@\tenrm}%
  \textfont1=\teni \scriptfont1=\seveni \scriptscriptfont1=\fivei
  \def\oldstyle{\fam\@ne\teni}%
  \textfont2=\tensy \scriptfont2=\sevensy \scriptscriptfont2=\fivesy
  \textfont\itfam=\tenit
  \def\it{\fam\itfam\tenit}%
  \textfont\slfam=\tensl
  \def\sl{\fam\slfam\tensl}%
  \textfont\bffam=\tenbf \scriptfont\bffam=\sevenbf
  \scriptscriptfont\bffam=\fivebf
  \def\bf{\fam\bffam\tenbf}%
  \textfont\ttfam=\tentt
  \def\tt{\fam\ttfam\tentt}%
  \abovedisplayskip=12pt plus 3pt minus 9pt
  \abovedisplayshortskip=0pt plus 3pt
  \belowdisplayskip=12pt plus 3pt minus 9pt
  \belowdisplayshortskip=7pt plus 3pt minus 4pt
  \smallskipamount=3pt plus 1pt minus 1pt
  \medskipamount=6pt plus 2pt minus 2pt
  \bigskipamount=12pt plus 4pt minus 4pt
  \normalbaselineskip=12pt
  \setbox\strutbox=\hbox{\vrule height8.5pt depth3.5pt width0pt}%
  \let\bigf@ntpc=\tenrm \let\smallf@ntpc=\sevenrm
  \let\petcap=\tenpc
  \normalbaselines\rm}

\def\eightpoint{%
  \textfont0=\eightrm \scriptfont0=\sixrm \scriptscriptfont0=\fiverm
  \def\rm{\fam\z@\eightrm}%
  \textfont1=\eighti \scriptfont1=\sixi \scriptscriptfont1=\fivei
  \def\oldstyle{\fam\@ne\eighti}%
  \textfont2=\eightsy \scriptfont2=\sixsy \scriptscriptfont2=\fivesy
  \textfont\itfam=\eightit
  \def\it{\fam\itfam\eightit}%
  \textfont\slfam=\eightsl
  \def\sl{\fam\slfam\eightsl}%
  \textfont\bffam=\eightbf \scriptfont\bffam=\sixbf
  \scriptscriptfont\bffam=\fivebf
  \def\bf{\fam\bffam\eightbf}%
  \textfont\ttfam=\eighttt
  \def\tt{\fam\ttfam\eighttt}%
  \abovedisplayskip=9pt plus 2pt minus 6pt
  \abovedisplayshortskip=0pt plus 2pt
  \belowdisplayskip=9pt plus 2pt minus 6pt
  \belowdisplayshortskip=5pt plus 2pt minus 3pt
  \smallskipamount=2pt plus 1pt minus 1pt
  \medskipamount=4pt plus 2pt minus 1pt
  \bigskipamount=9pt plus 3pt minus 3pt
  \normalbaselineskip=9pt
  \setbox\strutbox=\hbox{\vrule height7pt depth2pt width0pt}%
  \let\bigf@ntpc=\eightrm \let\smallf@ntpc=\sixrm
  \let\petcap=\eightpc
  \normalbaselines\rm}
\catcode`@=12

\tenpoint


\long\def\irmaaddress{{%
\bigskip
\eightpoint
\rightline{\quad
\vtop{\halign{\hfil##\hfil\cr
I.R.M.A. UMR 7501\cr
Universit\'e Louis Pasteur et CNRS,\cr
7, rue Ren\'e-Descartes\cr
F-67084 Strasbourg, France\cr
{\tt guoniu@math.u-strasbg.fr}\cr}}\quad}
}}



\catcode`\@=11
\def\pc#1#2|{{\bigf@ntpc #1\penalty \@MM\hskip\z@skip\smallf@ntpc%
	\uppercase{#2}}}
\catcode`\@=12

\def\pointir{\discretionary{.}{}{.\kern.35em---\kern.7em}\nobreak
   \hskip 0em plus .3em minus .4em }

\def\qed{\quad\raise -2pt\hbox{\vrule\vbox to 10pt{\hrule width 4pt
   \vfill\hrule}\vrule}}

\def\rem#1|{\par\medskip{{\it #1}\pointir}}

\def\vspace[#1]{\noalign{\vskip#1}}

\def\abstract#1{\vbox{\eightpoint\narrower\narrower 
\pc ABSTRACT|\pointir #1}}


\def\section#1{\goodbreak\par\vskip .3cm\centerline{\bf #1}
   \par\nobreak\vskip 3pt }

\long\def\th#1|#2\endth{\par\medbreak
   {\petcap #1\pointir}{\it #2}\par\medbreak}

\def\article#1|#2|#3|#4|#5|#6|#7|
    {{\leftskip=7mm\noindent
     \hangindent=2mm\hangafter=1
     \llap{{\tt [#1]}\hskip.35em}{\petcap#2}\pointir
     #3, {\sl #4}, {\bf #5} ({\oldstyle #6}),
     pp.\nobreak\ #7.\par}}
\def\livre#1|#2|#3|#4|
    {{\leftskip=7mm\noindent
    \hangindent=2mm\hangafter=1
    \llap{{\tt [#1]}\hskip.35em}{\petcap#2}\pointir
    {\sl #3}, #4.\par}}
\def\divers#1|#2|#3|
    {{\leftskip=7mm\noindent
    \hangindent=2mm\hangafter=1
     \llap{{\tt [#1]}\hskip.35em}{\petcap#2}\pointir
     #3.\par}}



\catcode`\@=11
\def\c@rr@#1{\vbox{%
  \hrule height \ep@isseur%
   \hbox{\vrule width\ep@isseur\vbox to \t@ille{%
           \vfil\hbox  to \t@ille{\hfil#1\hfil}\vfil}%
            \vrule width\ep@isseur}%
      \hrule height \ep@isseur}}
\def\ytableau#1#2#3#4{\vbox{%
  \gdef\ep@isseur{#2}
   \gdef\t@ille{#1}
    \def\\##1{\c@rr@{$#3 ##1$}}
  \lineskiplimit=-30cm \baselineskip=\t@ille%
    \advance \baselineskip by \ep@isseur%
     \halign{%
      \hfil$##$\hfil&&\kern -\ep@isseur%
       \hfil$##$\hfil \crcr#4\crcr}}}%
\catcode`\@=12

\def\Grille{\setbox13=\vbox to 5mm{\hrule width 110mm\vfill}
\setbox13=\vbox{\offinterlineskip
   \copy13\copy13\copy13\copy13\copy13\copy13\copy13\copy13
   \copy13\copy13\copy13\copy13\box13\hrule width 110mm}
\setbox14=\hbox to 5mm{\vrule height 65mm\hfill}
\setbox14=\hbox{\copy14\copy14\copy14\copy14\copy14\copy14
   \copy14\copy14\copy14\copy14\copy14\copy14\copy14\copy14
   \copy14\copy14\copy14\copy14\copy14\copy14\copy14\copy14\box14}
\ht14=0pt\dp14=0pt\wd14=0pt
\setbox13=\vbox to 0pt{\vss\box13\offinterlineskip\box14}
\wd13=0pt\box13}


\def\fleche(#1,#2)\dir(#3,#4)\long#5{%
\noalign{\nointerlineskip\leftput(#1,#2){\vector(#3,#4){#5}}\nointerlineskip}}


\def\hfl#1#2#3{\smash{\mathop{\hbox to#3{\rightarrowfill}}\limits
^{\scriptstyle#1}_{\scriptstyle#2}}}

\def\gfl#1#2#3{\smash{\mathop{\hbox to#3{\leftarrowfill}}\limits
^{\scriptstyle#1}_{\scriptstyle#2}}}


 \message{`lline' & `vector' macros from LaTeX}
 \catcode`@=11
\def\{{\relax\ifmmode\lbrace\else$\lbrace$\fi}
\def\}{\relax\ifmmode\rbrace\else$\rbrace$\fi}
\def\newcount{\alloc@0\count\countdef\insc@unt}
\def\newdimen{\alloc@1\dimen\dimendef\insc@unt}
\def\newwrite{\alloc@7\write\chardef\sixt@@n}

\newwrite\@unused
\newcount\@tempcnta
\newcount\@tempcntb
\newdimen\@tempdima
\newdimen\@tempdimb
\newbox\@tempboxa

\def\@spaces{\space\space\space\space}
\def\@whilenoop#1{}
\def\@whiledim#1\do #2{\ifdim #1\relax#2\@iwhiledim{#1\relax#2}\fi}
\def\@iwhiledim#1{\ifdim #1\let\@nextwhile=\@iwhiledim
        \else\let\@nextwhile=\@whilenoop\fi\@nextwhile{#1}}
\def\@badlinearg{\@latexerr{Bad \string\line\space or \string\vector
   \space argument}}
\def\@latexerr#1#2{\begingroup
\edef\@tempc{#2}\expandafter\errhelp\expandafter{\@tempc}%
\def\@eha{Your command was ignored.
^^JType \space I <command> <return> \space to replace it
  with another command,^^Jor \space <return> \space to continue without it.}
\def\@ehb{You've lost some text. \space \@ehc}
\def\@ehc{Try typing \space <return>
  \space to proceed.^^JIf that doesn't work, type \space X <return> \space to
  quit.}
\def\@ehd{You're in trouble here.  \space\@ehc}

\typeout{LaTeX error. \space See LaTeX manual for explanation.^^J
 \space\@spaces\@spaces\@spaces Type \space H <return> \space for
 immediate help.}\errmessage{#1}\endgroup}
\def\typeout#1{{\let\protect\string\immediate\write\@unused{#1}}}

\font\tenln    = line10
\font\tenlnw   = linew10

\newdimen\@wholewidth
\newdimen\@halfwidth
\newdimen\unitlength 

\unitlength =1pt


\def\thinlines{\let\@linefnt\tenln \let\@circlefnt\tencirc
  \@wholewidth\fontdimen8\tenln \@halfwidth .5\@wholewidth}
\def\thicklines{\let\@linefnt\tenlnw \let\@circlefnt\tencircw
  \@wholewidth\fontdimen8\tenlnw \@halfwidth .5\@wholewidth}

\def\linethickness#1{\@wholewidth #1\relax \@halfwidth .5\@wholewidth}

\newif\if@negarg

\def\lline(#1,#2)#3{\@xarg #1\relax \@yarg #2\relax
\@linelen=#3\unitlength
\ifnum\@xarg =0 \@vline
  \else \ifnum\@yarg =0 \@hline \else \@sline\fi
\fi}

\def\@sline{\ifnum\@xarg< 0 \@negargtrue \@xarg -\@xarg \@yyarg -\@yarg
  \else \@negargfalse \@yyarg \@yarg \fi
\ifnum \@yyarg >0 \@tempcnta\@yyarg \else \@tempcnta -\@yyarg \fi
\ifnum\@tempcnta>6 \@badlinearg\@tempcnta0 \fi
\setbox\@linechar\hbox{\@linefnt\@getlinechar(\@xarg,\@yyarg)}%
\ifnum \@yarg >0 \let\@upordown\raise \@clnht\z@
   \else\let\@upordown\lower \@clnht \ht\@linechar\fi
\@clnwd=\wd\@linechar
\if@negarg \hskip -\wd\@linechar \def\@tempa{\hskip -2\wd\@linechar}\else
     \let\@tempa\relax \fi
\@whiledim \@clnwd <\@linelen \do
  {\@upordown\@clnht\copy\@linechar
   \@tempa
   \advance\@clnht \ht\@linechar
   \advance\@clnwd \wd\@linechar}%
\advance\@clnht -\ht\@linechar
\advance\@clnwd -\wd\@linechar
\@tempdima\@linelen\advance\@tempdima -\@clnwd
\@tempdimb\@tempdima\advance\@tempdimb -\wd\@linechar
\if@negarg \hskip -\@tempdimb \else \hskip \@tempdimb \fi
\multiply\@tempdima \@m
\@tempcnta \@tempdima \@tempdima \wd\@linechar \divide\@tempcnta \@tempdima
\@tempdima \ht\@linechar \multiply\@tempdima \@tempcnta
\divide\@tempdima \@m
\advance\@clnht \@tempdima
\ifdim \@linelen <\wd\@linechar
   \hskip \wd\@linechar
  \else\@upordown\@clnht\copy\@linechar\fi}

\def\@hline{\ifnum \@xarg <0 \hskip -\@linelen \fi
\vrule height \@halfwidth depth \@halfwidth width \@linelen
\ifnum \@xarg <0 \hskip -\@linelen \fi}

\def\@getlinechar(#1,#2){\@tempcnta#1\relax\multiply\@tempcnta 8
\advance\@tempcnta -9 \ifnum #2>0 \advance\@tempcnta #2\relax\else
\advance\@tempcnta -#2\relax\advance\@tempcnta 64 \fi
\char\@tempcnta}

\def\vector(#1,#2)#3{\@xarg #1\relax \@yarg #2\relax
\@linelen=#3\unitlength
\ifnum\@xarg =0 \@vvector
  \else \ifnum\@yarg =0 \@hvector \else \@svector\fi
\fi}

\def\@hvector{\@hline\hbox to 0pt{\@linefnt
\ifnum \@xarg <0 \@getlarrow(1,0)\hss\else
    \hss\@getrarrow(1,0)\fi}}

\def\@vvector{\ifnum \@yarg <0 \@downvector \else \@upvector \fi}

\def\@svector{\@sline
\@tempcnta\@yarg \ifnum\@tempcnta <0 \@tempcnta=-\@tempcnta\fi
\ifnum\@tempcnta <5
  \hskip -\wd\@linechar
  \@upordown\@clnht \hbox{\@linefnt  \if@negarg
  \@getlarrow(\@xarg,\@yyarg) \else \@getrarrow(\@xarg,\@yyarg) \fi}%
\else\@badlinearg\fi}

\def\@getlarrow(#1,#2){\ifnum #2 =\z@ \@tempcnta='33\else
\@tempcnta=#1\relax\multiply\@tempcnta \sixt@@n \advance\@tempcnta
-9 \@tempcntb=#2\relax\multiply\@tempcntb \tw@
\ifnum \@tempcntb >0 \advance\@tempcnta \@tempcntb\relax
\else\advance\@tempcnta -\@tempcntb\advance\@tempcnta 64
\fi\fi\char\@tempcnta}

\def\@getrarrow(#1,#2){\@tempcntb=#2\relax
\ifnum\@tempcntb < 0 \@tempcntb=-\@tempcntb\relax\fi
\ifcase \@tempcntb\relax \@tempcnta='55 \or
\ifnum #1<3 \@tempcnta=#1\relax\multiply\@tempcnta
24 \advance\@tempcnta -6 \else \ifnum #1=3 \@tempcnta=49
\else\@tempcnta=58 \fi\fi\or
\ifnum #1<3 \@tempcnta=#1\relax\multiply\@tempcnta
24 \advance\@tempcnta -3 \else \@tempcnta=51\fi\or
\@tempcnta=#1\relax\multiply\@tempcnta
\sixt@@n \advance\@tempcnta -\tw@ \else
\@tempcnta=#1\relax\multiply\@tempcnta
\sixt@@n \advance\@tempcnta 7 \fi\ifnum #2<0 \advance\@tempcnta 64 \fi
\char\@tempcnta}

\def\@vline{\ifnum \@yarg <0 \@downline \else \@upline\fi}

\def\@upline{\hbox to \z@{\hskip -\@halfwidth \vrule
  width \@wholewidth height \@linelen depth \z@\hss}}

\def\@downline{\hbox to \z@{\hskip -\@halfwidth \vrule
  width \@wholewidth height \z@ depth \@linelen \hss}}

\def\@upvector{\@upline\setbox\@tempboxa\hbox{\@linefnt\char'66}\raise
     \@linelen \hbox to\z@{\lower \ht\@tempboxa\box\@tempboxa\hss}}

\def\@downvector{\@downline\lower \@linelen
      \hbox to \z@{\@linefnt\char'77\hss}}

\thinlines

\newcount\@xarg
\newcount\@yarg
\newcount\@yyarg
\newcount\@multicnt
\newdimen\@xdim
\newdimen\@ydim
\newbox\@linechar
\newdimen\@linelen
\newdimen\@clnwd
\newdimen\@clnht
\newdimen\@dashdim
\newbox\@dashbox
\newcount\@dashcnt
 \catcode`@=12


\newbox\tbox
\newbox\tboxa

\def\leftzer#1{\setbox\tbox=\hbox to 0pt{#1\hss}%
     \ht\tbox=0pt \dp\tbox=0pt \box\tbox}

\def\rightzer#1{\setbox\tbox=\hbox to 0pt{\hss #1}%
     \ht\tbox=0pt \dp\tbox=0pt \box\tbox}

\def\centerzer#1{\setbox\tbox=\hbox to 0pt{\hss #1\hss}%
     \ht\tbox=0pt \dp\tbox=0pt \box\tbox}

%
\def\image(#1,#2)#3{\vbox to #1{\offinterlineskip
    \vss #3 \vskip #2}}


\def\leftput(#1,#2)#3{\setbox\tboxa=\hbox{%
    \kern #1\unitlength
    \raise #2\unitlength\hbox{\leftzer{#3}}}%
    \ht\tboxa=0pt \wd\tboxa=0pt \dp\tboxa=0pt\box\tboxa}

\def\rightput(#1,#2)#3{\setbox\tboxa=\hbox{%
    \kern #1\unitlength
    \raise #2\unitlength\hbox{\rightzer{#3}}}%
    \ht\tboxa=0pt \wd\tboxa=0pt \dp\tboxa=0pt\box\tboxa}

\def\centerput(#1,#2)#3{\setbox\tboxa=\hbox{%
    \kern #1\unitlength
    \raise #2\unitlength\hbox{\centerzer{#3}}}%
    \ht\tboxa=0pt \wd\tboxa=0pt \dp\tboxa=0pt\box\tboxa}

\unitlength=1mm

\def\cput(#1,#2)#3{\noalign{\nointerlineskip\centerput(#1,#2){#3}
                             \nointerlineskip}}


\ifx\pdfoutput\jamaisdefined
\input epsf

\fi


\parskip 0pt plus 1pt

\def\article#1|#2|#3|#4|#5|#6|#7|
    {{\leftskip=7mm\noindent
     \hangindent=2mm\hangafter=1
     \llap{{\tt [#1]}\hskip.35em}{#2},\quad %
     #3, {\sl #4}, {\bf #5} ({\oldstyle #6}),
     pp.\nobreak\ #7.\par}}
\def\livre#1|#2|#3|#4|
    {{\leftskip=7mm\noindent
    \hangindent=2mm\hangafter=1
    \llap{{\tt [#1]}\hskip.35em}{#2},\quad %
    {\sl #3}, #4.\par}}
\def\divers#1|#2|#3|
    {{\leftskip=7mm\noindent
    \hangindent=2mm\hangafter=1
     \llap{{\tt [#1]}\hskip.35em}{#2},\quad %
     #3.\par}}


\def\runhookexp[#1]#2#3#4%
{\def\newline{\vskip 0pt\noindent \tt}
 \medskip
 \setbox1=\hbox{\kern#1{\vbox{#4}}\kern#1}%
  \dimen1=\ht1 \advance\dimen1 by #1 
  \dimen2=\dp1 \advance\dimen2 by #1 
  \setbox2=\hbox {\vbox to \dimen1{%
           \vbox{\hbox to #2{\vrule height 8pt\hfil #3\hfil}\hrule}\vfil}}
  \setbox1=\hbox{\vrule height\dimen1 depth\dimen2\box1\kern-#2\box2\vrule}%
  \setbox1=\vbox{\hrule\box1\hrule}%
  \advance\dimen1 by .4pt \ht1=\dimen1 
  \advance\dimen2 by .4pt \dp1=\dimen2  \box1\relax
}

\def\l{\lambda}

\def\mod{\mathop{\rm mod}}

\def\enslettre#1{\font\zzzz=msbm10 \hbox{\zzzz #1}}

\def\setN{\mathop{\enslettre N}}
\def\setP{\mathop{\cal P}}
\def\setH{\mathop{\cal H}}

\def\hmul{\mathop{\rm hmul}}

\def\setB{\mathop{\cal B}}
\def\setF{\mathop{\cal F}}
\def\setC{\mathop{\cal C}}

\def\secan{\mathop{\rm sec}}




\rightline{\Date}
\bigskip
\bigskip
\bigskip

\centerline{\bf Discovering hook length formulas by expansion technique}
\bigskip
\centerline{Guo-Niu HAN}
\bigskip\medskip

\abstract{
We introduce the hook length expansion technique and explain 
how to discover old and new hook length formulas for partitions and plane
trees. The new hook length formulas for trees obtained by our method can 
be proved rather easily, whereas those for partitions are much more difficult 
and 
some of them still remain open conjectures. We also develop a Maple package 
{\tt HookExp} for computing the hook length expansion.
The paper can be seen as a collection of hook length
formulas for partitons and plane trees.
All examples are illustrated by {\tt HookExp} and, for many easy cases,
expained by well-known combinatorial arguments.
}

\bigskip
\bigskip

\def\itemm#1|#2| {{\leftskip=7mm\noindent \hangindent=2mm\hangafter=1 
    \llap{\S#1.\ }{#2}\par}}

\centerline{\bf Summary}

\smallskip

\itemm1|Introduction. Selected hook formulas. Conjecture|
\itemm2|Classical hook length formulas for partitions.|
\itemm3|Hook length expansion algorithm and {\tt HookExp}.|
\itemm4|The exponent principle.|
\itemm5|Hook length formulas for partitions.|
\itemm6|Hook length formulas for binary trees.|
\itemm7|Hook length formulas for complete binary trees.|
\itemm8|Hook length formulas for Fibonacci trees.|


\def\sec{1}
\section{\sec. Introduction} 
The hook lengths for partitions and for plane trees play an important role
in Enumerative Combinatorics. The classical hook length formulas for those
two structures read 
$$
f_\l={n!\over \prod_{v\in \l} h_v}
\hbox{\quad and\quad}
f_T={n!\over \prod_{v\in T} h_v},
$$
where $f_\l$ (resp. $f_T$) is the number of standard Young tableaux
of shape $\l$ (resp. of increasing labeled binary trees of shape $T$).
See Sections 2 and~6 for notations and explanations.
From the above formulas we can derive 
$$
\sum_{\l\in\setP} x^{|\l|} \prod_{v\in\l} {1\over h_v^2}= e^{x}
\leqno{(\sec.1)}
$$
and
$$
\sum_{T\in\setB} x^{|\l|} \prod_{v\in T} {1\over h_v}= {1\over 1-x}.
\leqno{(\sec.2)}
$$
Formulas (\sec.1) and (\sec.2) are referred to as the basic 
{\it hook length formulas}, or {\it hook formulas}, for short.
\medskip

The numerous extensions or generalizations which have been proposed
in the literature led us to believe that a technical tool had to be
constructed that would make it possible to discover new hook
length formulas and also obtain the old ones in a systematic manner.
The purpose of this paper is to present such a tool that will
be called {\it hook length expansion} technique. In general, the new hook
length formulas for trees produced by that technique can be proved
easily, whereas those for partitions are much more difficult and 
some of them still remain open conjectures. 
\medskip

We also develop a Maple package {\tt HookExp}
for computing the hook length expansion, 
which can be downloaded freely from the author's web 
site.\footnote{$^{(*)}$}{Hook length formula homepage

{\tt http://math.u-strasbg.fr/\char126guoniu/hook}
} 
All the examples in the paper 
are illustrated by {\tt HookExp} and, for many easy cases,
explained by well-known combinatorial arguments.

\medskip

Sections 2-5 are devoted to the hook length formulas for partitions 
and Sections 6-8 for {\it plane trees}. 
Basic notions and classical hook length formulas for partitions 
are recalled in Section~2. Then, we introduce
the hook length expansion {\it algorithm} for partitions. 
In Section 4 we discuss some techniques
for discovering new hook length formulas, 
namely the {\it exponent principle}.
The new hook formulas 
for partitions $\l\in\setP$ 
(resp. for binary trees $T\in\setB$, 
for complete binary trees $T\in\setC$, 
for Fibonacci trees $T\in\setF$)
were suggested (but not proved!) by playing with the package
{\tt HookExp}. They are all collected in Section~5 (resp. Section 6, 7, 8).
The formulas we should like to single out are next stated. See
Sections 5-8 for notations, comments and/or proofs.

\proclaim Theorem \sec.1 [=5.5].
Let $t$ be a positive integer and 
$\hmul_t(\l)$ be the number of boxes $v$ such that $h_v(\l)$ is a multiple 
of $t$. Then
$$
\sum_{\l\in\setP} x^{|\l|} (-1)^{\hmul_t(\l)}
=
\prod_{k\geq 1}{
(1-x^{4tk})^t (1-x^{tk})^{2t}
\over
(1-x^{2tk})^{3t} (1-x^{k})
}.
$$

\proclaim Theorem \sec.2 [=5.9].
We have
$$
\sum_{\l\in \setP}x^{|\l|}\ \prod_{v\in\l}\ \bigl(1-{2 \over h_v^2}\bigr)
= 
\prod_{k\geq 1} { (1-x^k)}.
$$

\proclaim Theorem \sec.3 [=5.8].
We have
$$
\sum_{\l\in \setP}x^{|\l|} \prod_{v\in\l, h_v {\rm\, even}}\ 
\bigl(1-{2 \over h_v^2}\bigr)
=
\prod_{k\geq 1} { (1+x^k)}.
$$

The most general form of the above three theorems is Theorem 5.7. 
In fact, the latter theorem unifies several formulas, 
including the Jacobi triple 
product identity, the Macdonald identities for $A_\ell^{(a)}$, 
the generating functions for partitions (2.5) and for $t$-cores (5.4), 
the Nekrasov-Okounkov identity (5.1), Theorems 5.3-5.6 and Theorems~\sec.2-\sec.3.
See [Ha08e] for the proof and applications of Theorem 5.7.

\proclaim Conjecture \sec.4 [=5.2]. 
We have 
$$
\sum_{\l\in\setP} x^{|\l|} \prod_{v\in \l} \rho(z;h_v) 
= e^{x+zx^2/2},
$$
where the weight function $\rho(z; n)$ is defined by
$$
\rho(z; n)={\displaystyle
\sum_{k=0}^{\lfloor n/2\rfloor} {n\choose 2k} z^k \over
\displaystyle
n \sum_{k=0}^{\lfloor(n-1)/2\rfloor} {n\choose 2k+1} z^k 
}.
$$

\proclaim Theorem \sec.5 [=6.3].
We have
$$
\sum_{T\in\setB}x^{|T|}  \prod_{v\in T} {1\over h_v 2^{h_v-1}}
= e^x. 
$$

\proclaim Theorem \sec.6 [=6.6].
We have
$$
\sum_{T\in\setB}x^{|T|} \prod_{v\in T} {6\over h_v(h_v+2)}
= {1\over (1-x)^2}.
$$

\proclaim Theorem \sec.7 [=6.7].
We have
$$
\sum_{T\in\setB}x^{|T|} \prod_{v\in T} {h_v+3\over 2h_v}
= \bigl({1-\sqrt{1-4x}\over 2x}\bigr)^2. 
$$

The most general form of Theorem \sec.7 is Theorem 6.8.
In fact, the latter theorem unifies a lot of formulas, including 
the two classical hook formulas (6.4) and (6.5), Postnikov's formula (6.7)
and the generalization due to Lascoux, Du and Liu, another generalization
of Postnikov's formula (6.10), Theorems 6.6 and 6.7. 

\proclaim Theorem \sec.8 [=7.2].
We have
$$
\sum_{T\in\setC}x^{|T|}  \prod_{v\in T, h_v\geq 2} {1\over h_v2^{h_v-2}}= 
e^x.
$$

\proclaim Theorem \sec.9 [=8.6].
We have
$$
\sum_{T\in\setF}x^{|T|}  \prod_{v\in T, h_v\geq 2} 
{4(2h_v-1)(2h_v-3)\over (h_v+1)(5h_v-6)}
= 
{1-\sqrt{1-4x}\over 2x}.
$$
\medskip

For clarifying the nature of the paper we end the introduction by insisting 
on the following facts.

1. We introduce the hook length expansion technique by means of an explicit 
algorithm (Algorithm 3.1), together with the maple package {\tt HookExp}.

2. The package {\tt HookExp} is used to compute the first values
of the weight functions, which, in principle, suggest hook formulas to human
mathematicians. The package itself
does not output hook formulas, it does not prove either hook formulas!

3. We list new formulas found by {\tt HookEx}, 
but also some known formulas.
 
4. Most of the hook formulas for partitions are listed without any 
proofs. Instead,
we give the references containing the proofs, 
usually difficult and lengthy.

5. Most of the hook formulas for binary trees are listed with proofs, even
for well-known formulas,  because  most of them are proved in a unified way. 

6. Sometimes special cases of a master formula are also given, because
they 
have simpler forms with fewer parameters
and
show how the master formula was found by the author.

\def\sec{2}
\section{\sec. Classical hook length formulas for partitions} 
The basic notions needed here can be found in 
[Ma95, p.1; St99, p.287;  La01, p.1; Kn98, p.59; An76, p.1].
A {\it partition}~$\l$ is a sequence of positive 
integers $\l=(\l_1, \l_2,\cdots, \l_\ell)$ such that 
$\l_1\geq \l_2 \geq \cdots \geq \l_\ell>0$.
The integers
$(\l_i)_{i=1,2,\ldots, \ell}$ are called the {\it parts} of~$\l$,
the number $\ell$ of parts being the
{\it length} of $\l$ denoted by $\ell(\l)$.  
The sum of its parts $\l_1+ \l_2+\cdots+ \l_\ell$ is
denoted by $|\l|$.
Let $n$ be an integer, a partition 
$\l$ is said to be a partition of $n$ if $|\l|=n$. We write $\l\vdash n$.
The set of all partitions of $n$
is denoted by $\setP(n)$. 
The set of all partitions is denoted by~$\setP$,
so that $$\setP=\bigcup_{n\geq 0} \setP(n).$$
Each partition can be represented by its Ferrers diagram. For example,
$\l=(6,3,3,2)$ is a partition  and its Ferrers diagram is reproduced in 
Fig.~\sec.1.

{
\long\def\maplebegin#1\mapleend{}

\maplebegin

# --------------- begin maple ----------------------

# Copy the following text  to "makefig.mpl"
# then in maple > read("makefig.mpl");
# it will create a file "z_fig_by_maple.tex"

#\unitlength=1pt

Hu:= 12.4; # height quantities
Lu:= Hu; # large unity

X0:=-95.0; Y0:=15.6; # origin position

File:=fopen("z_fig_by_maple.tex", WRITE);

mhook:=proc(x,y,lenx, leny)
local i, d,sp, yy, xx, ct;
	sp:=Hu/8;
	ct:=0;
	for xx from x*Lu+X0 to x*Lu+X0+Hu by sp do
		yy := y*Hu+Y0; 
		fprintf(File, "\\vline(
				xx,   yy+sp*ct-0.2,    Lu*leny-sp*ct+0.1);
		ct:=ct+1;
	od:
	
	ct:=0;
	for yy from y*Hu+Y0 to y*Hu+Y0+Hu by sp do
		xx := x*Lu+X0; 
		fprintf(File, "\\hline(
				xx+sp*ct-0.2,   yy,    Lu*lenx-sp*ct+0.1);
		ct:=ct+1;
	od:

end;

mhook(1,1,2,3);

fclose(File);
# -------------------- end maple -------------------------
\mapleend

\setbox1=\hbox{$
\def\b{\\{\hbox{}}}
\ytableau{12pt}{0.4pt}{}
{\b &\b       \cr 
 \b &\b &\b   \cr
 \b &\b &\b   \cr
 \b &\b &\b &\b &\b &\b  \cr
\noalign{\vskip 3pt}
\noalign{\hbox{Fig.~\sec.1. Partition}}
}$
}
\setbox2=\hbox{$
\def\b{\\{\hbox{}}}
\unitlength=1pt%
\def\vline(#1,#2)#3|{\leftput(#1,#2){\lline(0,1){#3}}}%
\def\hline(#1,#2)#3|{\leftput(#1,#2){\lline(1,0){#3}}}%
\ytableau{12pt}{0.4pt}{}
{\b &\b       \cr 
 \b &\b &\b   \cr
 \b &\b &\b   \cr
 \b &\b &\b &\b &\b &\b  \cr
\noalign{\vskip 3pt}
\noalign{\hbox{Fig. \sec.2. Hook length}}%
}
\vline(-82.6,27.8)37.3|
\vline(-81.0,29.4)35.8|
\vline(-79.5,30.9)34.2|
\vline(-78.0,32.4)32.6|
\vline(-76.4,34.0)31.1|
\vline(-74.8,35.6)29.6|
\vline(-73.3,37.1)28.0|
\vline(-71.8,38.6)26.4|
\vline(-70.2,40.2)24.9|
\hline(-82.8,28.0)24.9|
\hline(-81.2,29.6)23.4|
\hline(-79.7,31.1)21.8|
\hline(-78.2,32.6)20.2|
\hline(-76.6,34.2)18.7|
\hline(-75.0,35.8)17.2|
\hline(-73.5,37.3)15.6|
\hline(-72.0,38.8)14.0|
\hline(-70.4,40.4)12.5|
$
}
\setbox3=\hbox{$
\ytableau{12pt}{0.4pt}{}
{\\2 &\\1       \cr 
 \\4 &\\3 &\\1   \cr
 \\5 &\\4 &\\2   \cr
 \\9 &\\8 &\\6 &\\3 &\\2 &\\1  \cr
\noalign{\vskip 3pt}
\noalign{\hbox{Fig. \sec.3. Hook lengths}}
}$
}
$$\box1\quad\box2\quad\box3$$
}

For each box $v$ in the Ferrers diagram of a partition $\l$, or
for each box $v$ in $\l$, for short, define the 
{\it hook length} of $v$, denoted by $h_v(\l)$ or $h_v$, to be the number of 
boxes $u$ such that  $u=v$,
or $u$ lies in the same column as $v$ and above $v$, or in the 
same row as $v$ and to the right of $v$ (see Fig.~\sec.2). 
The {\it hook length multi-set} of $\l$, denoted by $\setH(\l)$,
is the multi-set of all hook lengths of $\l$.
In Fig.~\sec.3 
the hook lengths of all boxes for the partition $\l=(6,3,3,2)$
have been written in each box. We have 
$\setH(\l)=\{2,1,4,3,1,5,4,2,9,8,6,3,2,1\}$.
The hook length plays an important role in Algebraic Combinatorics 
thanks to the famous hook formula
due to Frame, Robinson and Thrall [FRT54]
$$
f_\l={n!\over \prod_{h\in\setH(\l)} h}, \leqno{(\sec.1)}
$$
where $f_\l$ is the number of standard Young tableaux of shape $\l$
(see [St99, p.376; Kn98, p.59; GNW79; RW83; Ze84; GV85;  NPS97;  Kr99]).
\medskip

Recall that the Robinson-Schensted-Knuth 
correspondence (see, for example, [Kn98, p.49-59; St99, p.324])
is a bijection between the set of ordered pairs of standard Young tableaux 
of $\{1,2,\ldots, n\}$
of the same shape  and the set of 
permuations of order $n$. It provides a combinatorial proof of the following
identity.
$$
\sum_{\l\in n}\ f_\l^2 \ =\  n!  \leqno{(\sec.2)}
$$
By using (\sec.1) identity (\sec.2) can be written in the following
generating function form
$$
\sum_{\l\in\setP} x^{|\l|} \prod_{h\in \setH(\l)} {1\over h^2}
= e^x. \leqno{(\sec.3)}
$$
The Robinson-Schensted-Knuth 
correspondence also proves the fact that the number of
standard Young tableaux of $\{1,2,\ldots, \}$ 
is equal to the number of involutions  of order $n$ (see [Kn98b, p.47; Sch76]). 
In the generating function form this means that
$$
\sum_{\l\in\setP} x^{|\l|}\prod_{h\in \setH(\l)} {1\over h}
= e^{x+x^2/2}. \leqno{(\sec.4)} 
$$

The following identity is the well-known formula for the generating function 
of partitions [An76, p.3]. 
$$
\sum_{\l\in\setP} x^{|\l|} \prod_{h\in \setH(\l)} 1 
= \prod_{k\geq 1} {1\over  1-x^k}.  \leqno{(\sec.5)}
$$
In the present paper formulas (\sec.3), (\sec.4) and (\sec.5) are called also 
{\it hook formulas}. We will find other hook formulas in the next
sections.
\def\sec{3}
\section{\sec. Hook length expansion algorithm and {\tt HookExp}} 
For expressing our main algorithm in a handy manner it is convenient to 
introduce the following definition.
\medskip

{\it Definition \sec.1}. Let $\rho: \setN^* \rightarrow K$ be a map of the 
set of positive integers to some field $K$. 
Also let $f(x)\in K[[x]]$ be a formal power series in $x$ with 
coefficients in $K$ such that $f(0)=1$. If 
$$
\sum_{\l\in\setP} x^{|\l|} \prod_{h\in \setH(\l)} \rho(h) 
= f(x),  \leqno{(\sec.1)}
$$
the series $f(x)$ is called the {\it generating function} for partitions
by the weight function $\rho$. The left-hand side of (\sec.1) is called
the {\it hook length expansion} of $f(x)$. 
Furthermore, when both $\rho$ and $f(x)$ have simple (some people say ``nice")
form, 
equation (\sec.1) is called a {\it hook length formula}, or
{\it hook formula} for short.
\medskip

It is easy to see that the generating function $f(x)$ is uniquely determined
by the weight function $\rho$.  Conversely, the weight function $\rho$ 
can be uniquely determined by $f(x)$ in most cases. In the other cases (called
{\it singular cases}),
the weight function $\rho$ does not exist, or is not unique.
We next provide an algorithm for computing $\rho$ when $f(x)$ is given.
\medskip

Let $\setP_L(n)$ be the set of partitions $\l=(\l_1, \l_2, \ldots, \l_\ell)$ 
of $n$ such that $\ell(\l)=1$ or $\l_2=1$. The partitions in $\setP_L(n)$ are 
usually called {\it hooks}. The {\it hook length multi-set} $\setH(\l)$ of
a hook $\l$ of $n$ is simply 
$$
\setH(\l)=\{1,2,\cdots \ell(\l)-1,\  1,2,\cdots, n-\ell(\l),\  n\}.
\leqno{(\sec.2)}
$$
Let $\setP_Z(n)$ be the set of partitions
$\l=(\l_1, \l_2, \ldots, \l_\ell)$ of $n$ such that $\ell\geq 2$ and 
$\l_2\geq 2$. It is easy to see that 
the hook length multi-set of each partition of $\setP_Z(n)$
does not contain the integer $n$. 
Since $\setP(n)=\setP_L(n) \cup \setP_Z(n)$ we have
$$
\leqalignno{
\sum_{\l\vdash n} \prod_{h\in\setH(\l)} \rho(h)
&=
\sum_{\l\in \setP_L(n)} \prod_{h\in\setH(\l)} \rho(h)+
\sum_{\l\in \setP_Z(n)} \prod_{h\in\setH(\l)} \rho(h)\cr
&=\rho(n) \sum_{\l\in \setP_L(n)} \prod_{h=1}^{\ell(\l)-1} \rho(h) 
\prod_{h=1}^{n-\ell(\l)} \rho(h)+
\sum_{\l\in \setP_Z(n)} \prod_{h\in\setH(\l)} \rho(h). &{(\sec.3)}\cr
}
$$
The weight function $\rho$ can be obtained by the following 
algorithm.

\proclaim Algorithm \sec.1.
Let $f(x)=1+f_1 x +f_2 x^2 + f_3 x^3 + \cdots$ be a power series in $x$.
The weight function $\rho$ in the hook length expansion of $f(x)$ can be
calculated in the following manner. First, let $\rho(1)=f_1$. Then,
let $n\geq 2$ and suppose that
all values $\rho(k)$ for $1\leq k\leq n-1$ are known and
satisfy the following condition
$$D:=\sum_{\l\in \setP_L(n)} \prod_{h=1}^{\ell(\l)-1} \rho(h) 
\prod_{h=1}^{n-\ell(\l)} \rho(h) \not=0. \leqno{(\sec.4)}
$$
Then, by iteration, $\rho(n)$ is given by
$$ 
\rho(n)=
{f_n -\sum_{\l\in \setP_Z(n)} \prod_{h\in\setH(\l)} \rho(h) 
\over
D
}. \leqno{(\sec.5)}
$$

We only consider the power series $f(x)$ for which 
condition (\sec.4) holds. This is true in most cases. If for some reason
condition 
(\sec.4) fails to be true, 
we try to find an extension of $f(x)$ to avoid
the singularity.  More precisely, 
we try to find a series $F(x,t)\in K[[t]][[x]]$
such that $f(x)=F(x,0)$ and condition (\sec.4) holds for $F(x,t)$
(see (M.5.3) and (M.5.4) for an example). 

\medskip

The Maple package {\tt HookExp} is developed
for computing the first terms of the generating function $f(x)$ 
and the first values $\rho(n)$ in the hook length expansion. The underlying
variable of the series is always $x$. 
The input format for $f(x)$ is any
valid expression in Maple and the output format for $f(x)$ is
$$
1+f_1 x +f_2 x^2 + f_3 x^3 + f_4 x^ 4 + \cdots +f_n x^n.
$$
The input and output formats for $\rho(n)$ are the list
$$
[\rho(1), \rho(2), \rho(3), \ldots, \rho(n)].
$$
The procedure {\tt hookgen(rho)} computes the generating function $f(x)$ 
for the given weight function $\rho$, while the procedure {\tt hookexp(f, n)}
computes the weight function $\rho(k)$ for $k=1,2, \ldots, n$.
For example,
let us  verify identity (2.3) by using the {\tt HookExp} package.
\runhookexp[4pt]{18mm}{\hbox{\tt (M.\sec.1)}}{{
\newline > read("HookExp.mpl"): 
\newline > hooktype:="PA":   \# working on partitions
\newline > hookexp(exp(x), 8); 
$$
\bigl[1,{1 \over 4},{1 \over 9},{1 \over 16},{1 \over 25},
{1 \over 36},{1 \over 49},{1 \over 64}\bigr]
$$
\newline > hookgen(\%); 
$$
1+x+{1 \over 2}x^2+{1 \over 6}x^3+{1 \over 24}x^4+{1 \over 120}x^5+
{1 \over 720}x^6+{1 \over 5040}x^7+{1 \over 40320}x^8
$$
}} 
\medskip\noindent
Next,  verify identities (2.4) and (2.5).
\runhookexp[4pt]{18mm}{\hbox{\tt (M.\sec.2)}}{{
\newline > hookexp(exp(x+x\^{}2/2), 8); 
$$
\bigl[1,{1 \over 2},{1 \over 3},{1 \over 4},{1 \over 5},{1 \over 6},
{1 \over 7},{1 \over 8}\bigr]
$$
\newline > hookgen(\%); 
$$
1+x+x^2+{2 \over 3}x^3+{5 \over 12}x^4+{13 \over 60}x^5+{19 \over 180}x^6+
{29 \over 630}x^7+{191 \over 10080}x^8
$$
\newline > hookexp(product(1/(1-x\^{}k), k=1..9), 9); 
$$
[1,1,1,1,1,1,1,1,1]
$$
\newline > hookgen(\%); 
$$
1+x+2x^2+3x^3+5x^4+7x^5+11x^6+15x^7+22x^8+30x^9
$$
}} 
\medskip\noindent

\def\sec{4}
\section{\sec. The exponent principle} 
In principle, the 
{\tt HookExp} package gives rises to ``millions" of hook expansions.
But experience shows that only few of them can be dutifully named formulas. 
For example, with the very simple function $1/(1-x)$, we get the following
expansion.
\runhookexp[4pt]{18mm}{\hbox{\tt (M.\sec.1)}}{{
\newline > hookexp(1/(1-x), 8); 
$$
\bigl[1,{1 \over 2},{1 \over 2},{7 \over 12},{17 \over 25},
{447 \over 592},{160933 \over 197641},{105940688107 \over 124616941064}\bigr]
$$
}} 
\medskip \noindent
Apparently, no simple form can be obtained for $\rho(n)$.
Next, try to expand the generating function for the famous Catalan 
numbers (see, e.g., [St99, p.220]).
\runhookexp[4pt]{18mm}{\hbox{\tt (M.\sec.2)}}{{
\newline > hookexp((1-sqrt(1-4*x))/(2*x), 8); 
$$
\bigl[1,1,{5 \over 3},{37 \over 16},{823 \over 289},
{85028 \over 28605},{1055952653 \over 323028029}\bigr]
$$
}} 
\medskip \noindent
Not lucky again. 

Then, consider the generating function $f(x)$  for the given
weight function $\rho(n)=1+1/n$.
\runhookexp[4pt]{18mm}{\hbox{\tt (M.\sec.3)}}{{
\newline > hookgen([seq(1+1/n, n=1..8)]);
$$
1+2x+6x^2+{40 \over 3}x^3+31x^4+62x^5+{647 \over 5}x^6+
{3664 \over 15}x^7+{98467 \over 210}x^8
$$
}} 
\medskip \noindent
No evident formula for $f(x)$.
Those three examples tell us that it is not easy to discover
hook formula even with the help of {\tt HookExp}.
In fact, the author derived Algorithm 3.1 a long time ago, but 
never found any new hook length formula, until he recently discovered
the following {\it exponent principle}.

\proclaim The Exponent Principle.
If the power series $f(x)$ has a ``nice" hook length expansion, then there is
good chance that $f^z(x)$ has a ``nice" hook length expansion too.

The exponent principle was first discovered for binary trees. In such a case
the exponent principle can be partially justified (see (6.3)). 
It is then successfully applied for finding new hook length formulas for
partitions. The exponent principle for partitions has been verified
by experimental observation. However, the author
has no mathematical argument for proving or even partially explaning it. 

\medskip
Let us illustrate the exponent principle with the exponential function 
(see identity (2.3)).
\runhookexp[4pt]{18mm}{\hbox{\tt (M.\sec.4)}}{{
\newline > hookexp(exp(z*x), 8); 
$$
\bigl[z,{z \over 4},{z \over 9},{z \over 16},{z \over 25},
{z \over 36},{z \over 49},{z \over 64}\bigr]
$$
}} 
\medskip\noindent
It means that the following hook length expansion
$$
\sum_{\l\in\setP} x^{|\l|} \prod_{h\in \setH(\l)} {z\over h^2}
= e^{zx} \leqno{(\sec.1)}
$$
holds, but it is nothing new. We simply recover (2.3).
In the next section new hook length formulas
for partitions will be derived.
\def\sec{5}
\section{\sec. Hook length formulas for partitions} 
Let us apply the exponent principle to identity (2.5).
\runhookexp[4pt]{18mm}{\hbox{\tt (M.\sec.1)}}{{
\newline > hookexp(product(1/(1-x\^{}k)\^{}z, k=1..7), 7); 
$$
\bigl[z,{z+3 \over 4},{z+8 \over 9},{z+15 \over 16},
{z+24 \over 25},{z+35 \over 36},{z+48 \over 49}\bigr]
$$
}} 
\medskip\noindent
From the above expansion we derive the following hook length formula
for the power of Euler Product.
\proclaim Theorem \sec.1 [Nekrasov-Okounkov].
For any complex number $\beta$ we have 
$$
\sum_{\l\in \setP}\ \prod_{h\in\setH(\l)}\bigl(1-{\beta \over h^2}\bigr)x
\ =\ 
\prod_{k\geq 1} { (1-x^k)^{\beta-1}}.
\leqno{(\sec.1)}
$$

Theorem \sec.1 was discovered by Nekrasov and Okounkov in the study of 
the Seiberg-Witten Theory [NO06, {\it arXiv:hep-th/0306238v2}, 
formula (6.12), p.55].
In an unpublised paper (available on arXiv [Ha08a]) the author re-discovered
the Nekrasov-Okounkov identity (\sec.1) and gave an elementary proof by
using the Macdonald identities [Ma72]. 
Several applications ware also derived,
including the marked hook formula. 

\medskip

Now consider identity (2.4). The series
$$
f(x)=e^{x+x^2/2}
$$
is the generating function for the involutions.
By the exponent principle there is good chance that
$$
f^z(x)=e^{zx+zx^2/2}
\hbox{\quad or \quad }
f^{1/z}(zx)=e^{x+zx^2/2}
$$
has a ``nice" hook length expansion.
\runhookexp[4pt]{18mm}{\hbox{\tt (M.\sec.2)}}{{
\newline > hookexp(exp(x+z*x\^{}2/2), 9); 
$$
\leqalignno{
&\Bigl[1,\; {1+z \over 4},\; {3z+1 \over 9+3z},\; {z^2+6z+1 \over 16+16z},\;
   {5z^2+10z+1 \over 5z^2+50z+25},\;{z^3+15z^2+15z+1 \over 120z+36z^2+36},\cr
&\qquad\qquad\qquad{7z^3+35z^2+21z+1 \over 7z^3+147z^2+245z+49},\;
  {z^4+28z^3+70z^2+28z+1 \over 448z^2+64z^3+448z+64}\Bigr]\cr
}
$$
}} 
\medskip\noindent
The above values of $\rho$ suggests that the following new hook length formula, 
seen as an interpolation between permutations (2.3) 
and involutions (2.4), should hold.

\proclaim Conjecture \sec.2. 
We have 
$$
\sum_{\l\in\setP} x^{|\l|} \prod_{h\in \setH(\l)} \rho(z;h) 
= e^{x+zx^2/2}, \leqno{(\sec.2)}
$$
where the weight function $\rho(z; n)$ is defined by
$$
\rho(z; n)={\displaystyle
\sum_{k=0}^{\lfloor n/2\rfloor} {n\choose 2k} z^k \over
\displaystyle
n \sum_{k=0}^{\lfloor(n-1)/2\rfloor} {n\choose 2k+1} z^k 
}.
\leqno{(\sec.3)}
$$

When $z=1$, then $\rho(1; n)=1/n$. Identity (\sec.2) is true thanks to 
identity (2.4).  
When $z=0$, then $\rho(0; n)=1/n^2$. Identity (\sec.2) is also true since 
it becomes identity (2.3).  
However we cannot prove any other special cases of Conjecture \sec.2, 
except the above two values.  For more remarks about Conjecture \sec.2, 
see [Ha08d].

\medskip
Recall that a partition $\l$ is a $t$-core if the {\it hook length multi-set} 
of $\l$ does not contain the integer $t$.
It is known that the hook length multi-set of each $t$-core 
does not contain any {\it multiple} of $t$.
The generating function of the $t$-cores
is given by the following formula:
$$
\sum_\l x^{|\l|}
=\prod_{k\geq 1} {(1-x^{tk} )^t\over 1-x^k}, \leqno{(\sec.4)}
$$
where the sum ranges over all $t$-cores
[Kn98. p.69, p.612; St99, p.468; GKS90]. 
\medskip

Why do not expand the right-hand side of (\sec.4) by using
{\tt HookExp}? There is an interesting history hidden behind formula
(\sec.4). 
Since $t$ is a positive integer but not a free
parameter, we must choose a numerical value
for~$t$. Take $t=3$, formula (\sec.4) says that $\rho$ must 
have the following form
$$
[1,1,0,1,1,*,1,1,*,1,1,*],\leqno{(\sec.5)}
$$
where $*$ can be any numerical number.

\runhookexp[4pt]{18mm}{\hbox{\tt (M.\sec.3)}}{{ 
\newline > product((1-x\^{}(3*k))\^{}3/(1-x\^{}k), k=1..8):
\newline > hookexp(\%, 8); 
\newline {\tt Denominator is zero, no solution for n=8.}
\newline
\centerline{\tt [1, 1, 0, 1, 1, r[6], r[7], 0]}
}} 
\medskip\noindent
We cannot obtain (\sec.5) directly by using {\tt HookExp}. 
It is a sigular case. 
To avoid the sigularity
we replace the $3$ in the exponent by a free parameter $z$ 
(see comments after Algorithm 3.1).
Now {\tt hookexp} does not report
any error message. 
Finally we replace $z$ by $3$ to recover the $\rho$ shown in (\sec.5) 
(see also (M.\sec.12) for another variation of (M.\sec.3)).
\runhookexp[4pt]{18mm}{\hbox{\tt (M.\sec.4)}}{{ 
\newline > product((1-x\^{}(3*k))\^{}z/(1-x\^{}k), k=1..13):
\newline > hookexp(\%, 13); 
$$
\bigl[1,\ 1,\ 1-{z\over 3},\ 1,\ 1,\ 1-{z\over 12},\ 
1,\ 1,\ 1-{z\over 27},\ 1,\ 1,\ 1-{z\over 48},\ 1\bigr]
$$
\newline > subs(z=3, \%); 
$$\bigl[1,1,0,1,1,{3 \over 4},1,1,{8 \over 9},1,1,{15 \over 16},1\bigr]$$
}} 
\medskip\noindent
Moreover,
the above expansion suggests the following hook length formula,
which may be seen as an interpolation between identity (\sec.4) and
a specialization of (\sec.1):
$$
\sum_\l x^{|\l|} \prod_{v\in\l} \bigl(1- {t^2\over h_v^2}\bigr)
=\prod_{k\geq 1} {(1-x^k)^{t^2} \over 1-x^k}. \leqno{(\sec.6)}
$$

\proclaim Theorem \sec.3. 
We have the following hook length formula 
$$
\sum_{\l\in\setP} x^{|\l|} \prod_{h\in \setH(\l)} \rho(z;h) 
= \prod_{k\geq 1}{(1-x^{tk})^z\over 1-x^k} \leqno{(\sec.7)}
$$
where the weight function $\rho(z; n)$ is defined by
$$
\rho(z; n)=\cases{
1; &if $n\not\equiv 0 \mod t$.\cr
1-{tz\over n^2}; &if $n\equiv 0 \mod t$.\cr
}
\leqno{(\sec.8)}
$$

When $z=t$ we recover identity (\sec.4). When $t=1$ and $z=t^2$ we recover
identity (\sec.6).
Theorem \sec.3 is a special case of Theorem \sec.7, which is proved in [Ha08e].
\medskip

Now let us verify Theorem \sec.3 by using {\tt hookgen} 
for $t=2$ instead of {\tt hookexp}. In addition of
{\tt HookExp}, we also use the Maple package  {\tt qseries} developed by
Frank Garvan [Ga01]. Recall that the {\it Dedekind $\eta$-function}, is defined by 
$\eta(x)=x^{1/24} \prod_{m\geq 0} (1-x^m).$
 
\runhookexp[4pt]{18mm}{\hbox{\tt (M.\sec.5)}}{{ 
\newline > with(qseries);
\newline > r:=n-> if n mod 2=1 then 1 else 1-2*z/n\^{}2 fi: 
\newline > [seq(r(i), i=1..10)]; 
$$
\bigl[1,\  1-{z\over 2},\  1,\  1-{z\over 8},\  1,\  1-{z\over 18},\  1,\  1-{z\over 32},\  
1,\  1-{z\over 50},\  1\bigr]
$$
\newline > hookgen(\%): etamake(\%, x, 10): simplify(\%); 
$$
{x^{1/24-z/12}\ \eta(2\tau)^z \over \eta(\tau)}
$$
}} 
\medskip\noindent
As expected we obtain the right-hand side of (\sec.7) for $t=2$.
Next
we hope to obtain new hook formula by modifying slightly the 
above weight function.
Try to change the $1$ in odd position by~$-1$.
\runhookexp[4pt]{18mm}{\hbox{\tt (M.\sec.6)}}{{
\newline > r:=n-> if n mod 2=1 then -1 else 1-2*z/n\^{}2 fi: 
\newline > [seq(r(i), i=1..10)]; 
$$
\bigl[-1,\; 1-{z\over 2},\; -1,\; 1-{z\over 8},\; -1,\; 1-{z\over 18},\; -1,\; 1-{z\over 32},\; 
-1,\; 1-{z\over 50},\; -1 \bigr]
$$
\newline > hookgen(\%): etamake(\%, x, 10): simplify(\%); 
$$
x^{1/24-z/12}\ \eta(8\tau)^{2-z}\ \eta(4\tau)^{-5+3z}
\ \eta(2\tau)^{1-z}\ \eta(\tau)
$$
}} 
\medskip\noindent
The above expansion suggests the following hook length formula for partitions.
\proclaim Theorem \sec.4.
We have the following hook length formula 
$$
\sum_{\l\in\setP} x^{|\l|} \prod_{h\in \setH(\l)} \rho(z;h) 
= \prod_{k\geq 1}{
(1-x^k)(1-x^{4k})^{3z-5}
\over 
(1-x^{8k})^{z-2} (1-x^{2k})^{z-1}
}. \leqno{(\sec.9)}
$$
where the weight function $\rho(z; n)$ is defined by
$$
\rho(z; n)=\cases{
-1; &if $n\not\equiv 0 \mod 2$.\cr
1-{2z\over n^2}; &if $n\equiv 0 \mod 2$.\cr
}
\leqno{(\sec.10)}
$$

Inspired by Theorem \sec.4, we calculate the generating function for partitions
by the following periodical weight function~$\rho$.
\runhookexp[4pt]{18mm}{\hbox{\tt (M.\sec.7)}}{{
\newline > r:=n-> if n mod 3=0 then -1 else 1 fi: 
\newline > [seq(r(i), i=1..17)]; 
$$
[1, 1, -1, 1, 1, -1, 1, 1, -1, 1, 1, -1, 1, 1, -1, 1, 1]
$$
\newline > hookgen(\%): etamake(\%, x, 17): simplify(\%); 
$$
{x^{1/24}\eta(12\tau)^3\eta(3\tau)^6 \over \eta(6\tau)^9\eta(\tau)}
$$
}} 
\medskip\noindent
The above hook length expansion suggests the following formula. 
\proclaim Theorem \sec.5 [=1.1].
Let $t$ be a positive integer and 
$\hmul_t(\l)$ be the number of boxes $v$ such that $h_v(\l)$ is a multiple 
of $t$. Then
$$
\sum_{\l\in\setP} x^{|\l|} (-1)^{\hmul_t(\l)}
=
\prod_{k\geq 1}{
(1-x^{4tk})^t (1-x^{tk})^{2t}
\over
(1-x^{2tk})^{3t} (1-x^{k})
}. \leqno{(\sec.11)}
$$

In fact, Theorem \sec.5 can be generalized by replacing $-1$ by $z$. 
\runhookexp[4pt]{18mm}{\hbox{\tt (M.\sec.8)}}{{ 
\newline > f := k -> (1-x\^{}(3*k))\^{}3/(1-(z*x\^{}3)\^{}k)\^{}3/(1-x\^{}k):
\newline > hookexp(product(f(k),k=1..15), 15); 
$$
[1, 1, z, 1, 1, z, 1, 1, z, 1, 1, z, 1, 1, z]
$$
}} 
\medskip\noindent
\proclaim Theorem \sec.6.
Let $t$ be a positive integer. Then
$$
\sum_{\l\in\setP} x^{|\l|} z^{\hmul_t(\l)}
=
\prod_{k\geq 1}{
(1-x^{tk})^{t}
\over
(1-(zx^t)^k)^{t} (1-x^{k})
}. \leqno{(\sec.12)}
$$

We can unify (M.\sec.4) and (M.\sec.8) in the following manner.
\runhookexp[4pt]{18mm}{\hbox{\tt (M.\sec.9)}}{{ 
\newline > N:=14: t:=3:
\newline > r:=n-> if n mod t=0 then  y*(1-t*z/n\^{}2) else 1 fi:
\newline > fk := k-> %
  ((1-x\^{}(t*k))\^{}t)/((1-(y*x\^{}t)\^{}k)\^{}(t-z))/(1-x\^{}k):  
\newline > [seq(r(i), i=1..N)];
$$
\bigl[1,\; 1,\; y-{yz \over 3},\; 1,\; 1,\; y-{yz \over 12},\; 1,\; 1,\; 
y-{yz \over 27},\; 1,\; 1,\; y-{yz \over 48},\; 1,\; 1\bigr]
$$
\newline > hookgen(\%) - product(fk(k), k=1..N):
\newline > series(\%,x,N+1): simplify(\%);
$$
O(x^{15})
$$
}} 
\medskip\noindent
The above expansion suggests the following hook length formula.
\proclaim Theorem \sec.7. 
We have 
$$
\sum_{\l\in\setP} x^{|\l|} \prod_{h\in \setH(\l)} \rho(z;h) 
= 
\prod_{k\geq 1}
{
(1-x^{tk})^t
\over
(1-(yx^t)^k)^{t-z}(1-x^k)  
},\leqno{(\sec.13)}
$$
where the weight function $\rho(z; n)$ is defined by
$$
\rho(z; n)=\cases{
1; &if $n\not\equiv 0 \mod t$.\cr
y-{tyz\over n^2}; &if $n\equiv 0 \mod t$.\cr
}
\leqno{(\sec.14)}
$$

The proof of Theorem \sec.7, as well as some applications can be found  
in [Ha08e].
Let us single out the very simple case when $t=2, y=z=1$. 
\runhookexp[4pt]{18mm}{\hbox{\tt (M.\sec.10)}}{{ 
\newline > hookexp(product(1+x\^{}k, k=1..14),14);
$$
\bigl[1,{1 \over 2},1,{7 \over 8},1,{17 \over 18},1,{31 \over 32},1,
{49 \over 50},1,{71 \over 72},1,{97 \over 98}\bigr]
$$
}} 
\proclaim Theorem \sec.8 [=1.3]. 
We have the following hook length formula 
$$
\sum_{\l\in\setP} x^{|\l|} \prod_{h\in \setH(\l), h {\rm\ even}} 
\bigl(1-{2\over h^2}   \bigr) 
= 
\prod_{k\geq 1} (1+x^k). \leqno{(\sec.15)}
$$

The above theorem is to be compared with the following specialization of
Theorem \sec.7 when $z=2, y=t=1$.
\proclaim Theorem \sec.9 [=1.2]. 
We have the following hook length formula 
$$
\sum_{\l\in\setP} x^{|\l|} \prod_{h\in \setH(\l)} 
\bigl(1-{2\over h^2}   \bigr) 
= 
\prod_{k\geq 1} (1-x^k). \leqno{(\sec.16)}
$$

There are also other hook formulas that are not specialization of 
Theorem~\sec.7.
Consider the weight function $\rho$ that counts the corners 
(their hook lengths are $1$) of partitions. 
\runhookexp[4pt]{18mm}{\hbox{\tt (M.\sec.11)}}{{ 
\newline > [z,seq(1, i=1..7)]; 
$$
[z, 1, 1, 1, 1, 1, 1, 1]
$$
\newline > hookgen(\%);
$$
\leqalignno{
&1+(z)x+(2z)x^2+(2z+z^2)x^3+(3z+2z^2)x^4+(2z+5z^2)x^5\cr
&+(4z+6z^2+z^3)x^6+(2z+11z^2+2z^3)x^7+(4z+13z^2+5z^3)x^8\cr
}
$$
}} 
\medskip\noindent
The above generating function corresponds to the
sequence A116608 in the on-line encyclopedia of integer sequences [Slo]
and is equal to the right-hand side of (\sec.17) below.
\proclaim Theorem \sec.10. 
We have 
$$
\sum_{\l\in\setP} x^{|\l|} \prod_{h\in \setH(\l), h =1}  z
= 
\prod_{k\geq 1}{1+(z-1)x^k\over 1-x^k}.
\leqno{(\sec.17)}
$$

On the other hand, 
take (M.\sec.3) and change the ``$-$" in the numerator by ``$+$", we get 
\runhookexp[4pt]{18mm}{\hbox{\tt (M.\sec.12)}}{{ 
\newline > product((1+x\^{}(3*k))\^{}3/(1-x\^{}k), k=1..8):
\newline > hookexp(\%, 8); 
\newline
\centerline{\tt [1, 1, 2, 1, 1, 1, 1, 1]}
}} 
\medskip\noindent
The above expansion suggests the following formula.
\proclaim Theorem \sec.11. 
We have 
$$
\sum_{\l\in\setP} x^{|\l|} \prod_{h\in \setH(\l), h =t}  2
= 
\prod_{k\geq 1}
{
(1+x^{tk})^t
\over
1-x^k
}.\leqno{(\sec.18)}
$$

In [Ha08e] we proved the following unified form of Theorems \sec.10 and \sec.11
by using
the properties of a classical bijection
which maps each partition to its $t$-core and $t$-quotient
[Ma95, p.12; St99, p.468; JK81, p.75; GSK90]. 

\proclaim Theorem \sec.12. 
For any complex number $z$ we have 
$$
\sum_{\l\in\setP} x^{|\l|} \prod_{h\in \setH(\l), h =t}  z
= 
\prod_{k\geq 1}
{
(1+(z-1)x^{tk})^t
\over
1-x^k
}.\leqno{(\sec.19)}
$$

\def\sec{6}
\section{\sec. Hook length formulas for binary trees} 

The basic notions for binary trees can be found in 
[St97,p.295; Kn98a, p.308-313; Vi81]. A {\it binary tree} $T$ with $n$ vertices
is defined {\it recursively} as follows. 
Either $T$ is empty, or else one specially designated 
vertex $v$ is called the root of $T$,  and the remaining vertices 
(excluding the root) are put into an ordered pair $(T', T'')$
of binary trees (possibly empty), which are called {\it subtrees}
of the root $v$. 
The {\it hook length} of the root $v$, denoted by $h_v(T)$ or $h_v$, is 
just the number of vertices $n$, which is also denoted by $|T|$.
Each vertex is called {\it leaf} if his two subtrees
are both empty. 
The {\it hook length multi-set}
$\setH(T)=\{h_v\mid v\in T\}$
of $T$
is defined to be the multi-set of hook lengths of all vertices $v$ of $T$. 
Finally,  let $\setB$ (resp. $\setB(n)$) denote
the set of all binary trees (resp. all binary trees with $n$ vertices),
so that 
$$\setB=\bigcup_{n\geq 0} \setB(n).$$
For example, there are five binary trees with $n=3$ vertices.
%

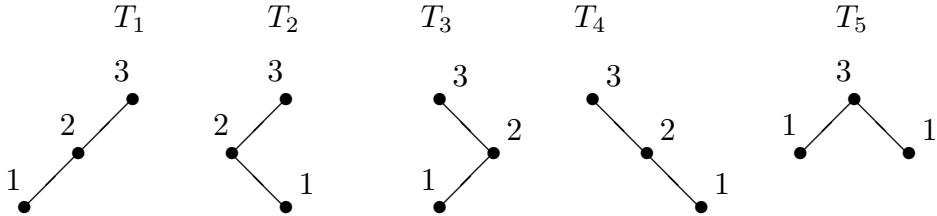
\long\def\maplebegin#1\mapleend{}

\maplebegin

# --------------- begin maple ----------------------

# Copy the following text  to "makefig.mpl"
# then in maple > read("makefig.mpl");
# it will create a file "z_fig_by_maple.tex"

#\unitlength=1pt

Hu:= 6; # height quantities

X0:=5.0; Y0:=0; # origin position

File:=fopen("z_fig_by_maple.tex", WRITE);

pline:=proc(x,y) # X0,Y0 = offset
local a,b,len;
	len:=1;
	fprintf(File, "\\pline(
end;

nline:=proc(x,y) # X0,Y0 = offset
local a,b,len;
	len:=1;
	fprintf(File, "\\nline(
end;

mydot:=proc(x,y) # X0,Y0 = offset
local a,b,len;
	len:=1;
	fprintf(File, "\\mydot(
end;

mylabel:=proc(x,y, text) # X0,Y0 = offset
local a,b,len;
	len:=1;
	fprintf(File, "\\mylabel(
end;

mylabel2:=proc(x,y, text) # X0,Y0 = offset
local a,b,len;
	len:=1;
	fprintf(File, "\\mylabel(
end;

dotlabel:=proc(x,y,t) mydot(x,y); mylabel(x,y,t); end;
dotlabel2:=proc(x,y,t) mydot(x,y); mylabel2(x,y,t); end;

DXX:=23;

X0:=-2.6*DXX;
pline(2,1); pline(3,2); 
dotlabel(2,1, "1"); dotlabel(3,2, "2"); dotlabel(4,3, "3");
mylabel(4,4, "$T_1$");

X0:=X0+DXX;

nline(2,2); pline(2,2); 
dotlabel2(3,1, "1"); dotlabel(2,2, "2"); dotlabel(3,3, "3");
mylabel(3,4, "$T_2$");

X0:=X0+DXX;
nline(2,3); pline(2,1); 
dotlabel(2,1, "1"); dotlabel2(3,2, "2"); dotlabel2(2,3, "3");
mylabel(2,4, "$T_3$");

X0:=X0+DXX; 
nline(1,3); nline(2,2); 
dotlabel2(3,1, "1"); dotlabel2(2,2, "2"); dotlabel2(1,3, "3");
mylabel(1,4, "$T_4$");

X0:=X0+DXX;
pline(1,2); nline(2,3); 
dotlabel(1,2, "1"); dotlabel2(3,2, "1"); dotlabel(2,3, "3");
mylabel(2,4, "$T_5$");

fclose(File);

# -------------------- end maple -------------------------

\mapleend


\newbox\boxarbre
\def\pline(#1,#2)#3|{\leftput(#1,#2){\lline(1,1){#3}}}
\def\nline(#1,#2)#3|{\leftput(#1,#2){\lline(1,-1){#3}}}
\def\mydot(#1,#2)|{\leftput(#1,#2){$\bullet$}}
\def\mylabel(#1,#2)#3|{\leftput(#1,#2){#3}}
\setbox\boxarbre=\vbox{\vskip
30mm\offinterlineskip 
%
\pline(-47.8,6.0)6.0|
\pline(-41.8,12.0)6.0|
\mydot(-48.6,5.2)|
\mylabel(-49.8,8.0){1}|
\mydot(-42.6,11.2)|
\mylabel(-43.8,14.0){2}|
\mydot(-36.6,17.2)|
\mylabel(-37.8,20.0){3}|
\mylabel(-37.8,26.0){$T_1$}|
\nline(-24.8,12.0)6.0|
\pline(-24.8,12.0)6.0|
\mydot(-19.6,5.2)|
\mylabel(-17.3,7.5){1}|
\mydot(-25.6,11.2)|
\mylabel(-26.8,14.0){2}|
\mydot(-19.6,17.2)|
\mylabel(-20.8,20.0){3}|
\mylabel(-20.8,26.0){$T_2$}|
\nline(-1.8,18.0)6.0|
\pline(-1.8,6.0)6.0|
\mydot(-2.6,5.2)|
\mylabel(-3.8,8.0){1}|
\mydot(3.4,11.2)|
\mylabel(5.7,13.5){2}|
\mydot(-2.6,17.2)|
\mylabel(-0.3,19.5){3}|
\mylabel(-3.8,26.0){$T_3$}|
\nline(15.2,18.0)6.0|
\nline(21.2,12.0)6.0|
\mydot(26.4,5.2)|
\mylabel(28.7,7.5){1}|
\mydot(20.4,11.2)|
\mylabel(22.7,13.5){2}|
\mydot(14.4,17.2)|
\mylabel(16.7,19.5){3}|
\mylabel(13.2,26.0){$T_4$}|
\pline(38.2,12.0)6.0|
\nline(44.2,18.0)6.0|
\mydot(37.4,11.2)|
\mylabel(36.2,14.0){1}|
\mydot(49.4,11.2)|
\mylabel(51.7,13.5){1}|
\mydot(43.4,17.2)|
\mylabel(42.2,20.0){3}|
\mylabel(42.2,26.0){$T_5$}|%
}
$$
\kern-4mm\box\boxarbre
$$
We have $\setH(T_1)=\setH(T_2)= \setH(T_3)= \setH(T_4)= \{1,2,3\}$
and $\setH(T_5)=\{1,1,3\}$.

As done for the partitions in Definition 3.1 we define the {\it hook
length expansion} for binary trees by
$$
\sum_{T\in\setB} x^{|T|} \prod_{h\in \setH(T)} \rho(h) 
= f(x),  \leqno{(\sec.1)}
$$
where $f(x)\in K[[x]]$ is a power series in $x$ with coefficients in $K$ such 
that $f(0)=1$.
See Section 3 for other related definitions and comments about hook length
expansion. For computing the weight function $\rho$, we need find an analogue 
of Algorithm 3.1 for binary trees. No surprise, it is much easier to find
a formula for computing $\rho$, since the binary tree structure is more simple 
compared with the partition structure.

\medskip
Let $f(x)=1+f_1x +f_2x^2+f_3x^3+\cdots$ be the generating function for 
binary trees by the weight function  $\rho$.
With each $T\in\setB(n)$ ($n\geq 1$) we can 
associate a triplet $(T', T'', v)$, where
$T'\in\setB(k)$ ($0\leq k\leq n-1$),
$T''\in\setB(n-1-k)$  and the root $v$ of $T$ whose hook length $h_v=n$.
Hence (\sec.1) is equivalent to
$$
\rho(n)\sum_{k=0}^{n-1} f_k f_{n-1-k}=f_n \quad (n\geq 1).
\leqno{(\sec.2)}
$$
Formula (\sec.2) can be used to calculate $f(x)$ for a given $\rho$,
or to calculate $\rho$ for a given $f(x)$. It has also
the equivalent form
$$
\rho(n) = {[x^n] f(x) \over [x^{n-1}] f^2(x)}, \leqno{(\sec.3)}
$$
where $[x^n]f(x)$ means the coefficient of $x^n$ in the power series $f(x)$.
From (\sec.3) we may say that finding a {\it hook length formula} is 
equivalent to finding a formal power
series $f(x)$ such that $[x^n] f(x)/[x^{n-1}] f^2(x)$ has a 
``nice" form in $n$. 
\medskip
Next we use the maple package {\tt HookExp} to find hook formulas for 
binary trees. The syntax of the two procedures {\tt hookexp} and
{\tt hookgen} are the same as for partitions. The rest of this section
contains some sessions, and 
each session contains
three parts: 
(i)~experiment with {\tt HookExp}; 
(ii)~hook formula suggested by the experiment; 
(iii)~proof and/or comments of the hook formula. 
All proofs of the hook formulas presented in this section
are always based on relation (\sec.3). 

\runhookexp[4pt]{18mm}{\hbox{\tt (M.\sec.1)}}{{
\newline > hooktype:="BT":  \# working on binary trees
\newline > hookexp(1/(1-x), 9); 
$$
\bigl[1,{1 \over 2},{1 \over 3},{1 \over 4},{1 \over 5},{1 \over 6},
{1 \over 7},{1 \over 8},{1 \over 9}\bigr]
$$
}} 
\medskip\noindent
\proclaim Theorem \sec.1.
We have
$$
\sum_{T\in\setB}x^{|T|}  \prod_{h\in \setH(T)} {1\over h}= 
{1\over 1-x}. \leqno{(\sec.4)}
$$

{\it Proof}. From (\sec.3)
$$
\rho(n) = {[x^n] 1/(1-x) \over [x^{n-1}] 1/(1-x)^2}=1/n.\qed
$$

{\it Remark}.
It is well-known [Kn98b,p.67; St75] that
the number of ways to label the vertices of $T$ with $\{1,2,\ldots, n\}$,
such that the label of each vertex is less than that of its descendants
(called {\it increasing labeled binary trees}, or {\it labeled binary trees} 
for short),
is equal to $n!$ divided by the product of the $h_v$'s ($v\in T$).
On the other hand, each labeled binary tree with $n$ vertices
is in bijection with a permutation of order $n$ [St97,p.24;FS73;Vi81], 
so that
$$
\sum_{T\in\setB(n)} n! \prod_{v\in T} {1\over h_v}= n! 
$$
This gives a combinatorial proof of Theorem \sec.1.

\runhookexp[4pt]{18mm}{\hbox{\tt (M.\sec.2)}}{{
\newline > hookexp((1-sqrt(1-4*x))/(2*x), 9); 
$$
[1, 1, 1, 1, 1, 1, 1, 1]
$$
}} 
\medskip\noindent
\proclaim Theorem \sec.2.
We have
$$
\sum_{T\in\setB}x^{|T|}  \prod_{h\in \setH(T)} {1}
= {1-\sqrt{1-4x}\over 2x}. 
\leqno{(\sec.5)}
$$

{\it Proof}.
Let $f(x)$ be the right-hand side of (\sec.5). Then
$$
[x^{n-1}] f^2(x) =
[x^{n-1}] (f(x)-1)/x = [x^n] f(x). \qed
$$

{\it Remark}.
Formula (\sec.5) implies 
that the number of binary trees with $n$ vertices is equal to the
$n$-th Catalan number (see, e.g., [St99, p.220]) 
$$
\sum_{T\in\setB(n)} 1= {1\over n+1} {2n\choose n}. 
$$

In the following experiments 
we make use of the maple package {\tt Guess}, translated by 
B\'eraud and Gauthier [BG04] from the Mathematica package {\tt Rate}
devoloped by Krattenthaler [Kr01]. 
\runhookexp[4pt]{18mm}{\hbox{\tt (M.\sec.3)}}{{
\newline > with(GUESS): 
\newline > guess:=proc(r) subs(\{\_i[0]=n, \_i[1]=k, \_i[2]=m\}, Guess(r)):
\newline \hskip 31mm   simplify(\%); op(\%); end:
\newline > hookexp(exp(x), 9);  
$$
\bigl[1,{1 \over 4},{1 \over 12},{1 \over 32},{1 \over 80},
{1 \over 192},{1 \over 448},{1 \over 1024},{1 \over 2304}\bigr]
$$
\newline > guess(\%);
$$
{{2^{1-n}} \over n}
$$
}} 
\medskip\noindent
\proclaim Theorem \sec.3.
We have
$$
\sum_{T\in\setB}x^{|T|}  \prod_{h\in \setH(T)} {1\over h 2^{h-1}}
= e^x. 
\leqno{(\sec.6)}
$$

{\it Proof}.
By (\sec.3)
$$
\rho(n)={[x^n] e^x \over [x^{n-1}] e^{2x}} = {1/n! \over 2^{n-1}/(n-1)!} 
={1\over n2^{n-1}}.\qed
$$
{\it Remark}. We do not have any combinatorial proof of Theorem \sec.3. 
See also [Ha08b].\qed

\runhookexp[4pt]{18mm}{\hbox{\tt (M.\sec.4)}}{{
\newline > [seq(1+1/n, n=1..7)];  
$$
\bigl[2,{3 \over 2},{4 \over 3},{5 \over 4},{6 \over 5},{7 \over 6},
{8 \over 7}\bigr]
$$
\newline > f:=hookgen(\%);
$$
1+2x+6x^2+{64 \over 3}x^3+{250 \over 3}x^4+{1728 \over 5}x^5+
{67228 \over 45}x^6+{2097152 \over 315}x^7
$$
}} 
\medskip\noindent
\proclaim Theorem \sec.4 [Postnikov].
We have
$$
\sum_{T\in\setB} x^{|T|}
\prod_{h\in \setH(T)} \bigl(1+{1\over h}\bigr)
= \sum_{n\geq 0}(n+1)^{n-1} {(2x)^n\over n!}.
\leqno{(\sec.7)}
$$

{\it Proof}.
Let $G(x)$ be a power series such that
$$
G(x)=\exp(x G(x)). \leqno{(\sec.8)}
$$
By the Lagrange inversion formula $G(x)^z$ has the following explicit 
expansion:
$$
G(x)^{z}=\sum_{n\geq 0} z(n+z)^{n-1} {x^n\over n!}. \leqno{(\sec.9)}
$$
The right-hand side of (\sec.7) is $G(2x)$. By (\sec.3)
$$
\rho(n)={[x^n] G(2x) \over [x^{n-1}] G^2(2x)} = 
{(n+1)^{n-1} 2^n/n! \over 2(n+1)^{n-2} 2^{n-1}/(n-1)!} 
=1+{1\over n}.
$$
Further combinatorial proofs and extensions 
have been proposed by several authors 
[Po04, CY08, DL08, GS06, MY07, Se08, Ha08c].~\qed

\runhookexp[4pt]{18mm}{\hbox{\tt (M.\sec.5)}}{{
\newline > [seq(1+1/n, n=1..9)];  
$$
\bigl[2,{3 \over 2},{4 \over 3},{5 \over 4},{6 \over 5},
{7 \over 6},{8 \over 7},{9 \over 8},{10 \over 9}\bigr]
$$
\newline > f:=hookgen(\%): hookexp(f\^{}z, 7): map(factor, \%);
$$
\bigl[2z,{2+z \over 2},{(z+3)^2 \over 6z+6},{(z+4)^3 \over 4(2z+3)^2},
{(z+5)^4 \over 40(2+z)^3},{(z+6)^5 \over 6(2z+5)^4},
{(z+7)^6 \over 224(z+3)^5}\bigr]
$$
}} 
\medskip\noindent
\proclaim Theorem \sec.5.
We have
$$
\sum_{T\in\setB}x^{|T|} \prod_{v\in T} {(z+h)^{h-1}\over h (2z+h-1)^{h-2}}
= \sum_{n\geq 0}{z}(z+n)^{n-1} {(2x)^n\over n!}.
\leqno{(\sec.10)}
$$

{\it Proof}.
From (\sec.3), (\sec.8) and (\sec.9) we have
$$
\leqalignno{
\rho(n)
&={[x^n] G^z(2x) \over [x^{n-1}] G^{2z}(2x)} 
= 
{z(n+z)^{n-1} 2^n/n! \over 2z(n+2z-1)^{n-2} 2^{n-1}/(n-1)!} \cr
&= {(z+n)^{n-1}\over n (2z+n-1)^{n-2}}.\qed\cr
}
$$
{\it Remark}. We do not have any combinatorial proof of Theorem \sec.5. 
See also [Ha08c].

\runhookexp[4pt]{18mm}{\hbox{\tt (M.\sec.6)}}{{
\newline > hookexp(1/(1-x)\^{}2, 9);  
$$
\bigl[2,{3 \over 4},{2 \over 5},{1 \over 4},{6 \over 35},
{1 \over 8},{2 \over 21},{3 \over 40},{2 \over 33}\bigr]
$$
\newline > guess(\%);
$$
{6 \over n(n+2)}
$$
}} 
\medskip\noindent
\proclaim Theorem \sec.6.
We have
$$
\sum_{T\in\setB}x^{|T|} \prod_{v\in T} {6\over n(n+2)}
= {1\over (1-x)^2}.
\leqno{(\sec.11)}
$$
More generally,
$$
\sum_{T\in\setB}x^{|T|} \prod_{h\in \setH(T)} {
\prod_{i=1}^{h-1} (z+i) \over
2h \prod_{i=1}^{h-2} (2z+i)     }
= {1\over (1-x)^z}
\leqno{(\sec.12)}
$$
or
$$
\sum_{T\in\setB(n)} \prod_{h\in \setH(T)} {
\prod_{i=1}^{h-1} (z+i) \over
2h \prod_{i=1}^{h-2} (2z+i)     }
= {1\over n!} \prod_{i=0}^{n-1}(z+i).
$$

{\it Proof}.
By (\sec.3)
$$
\rho(n)
={[x^n] 1/(1-x)^z \over [x^{n-1}] 1/(1-x)^{2z}} 
={{n+z-1 \choose n} \over {n+2z-2 \choose n-1}}.\qed
$$

\runhookexp[4pt]{18mm}{\hbox{\tt (M.\sec.7)}}{{
\newline > hookexp(((1-sqrt(1-4*x))/(2*x))\^{}2, 10); 
$$
\bigl[2,{5 \over 4},1,{7 \over 8},{4 \over 5},{3 \over 4},
{5 \over 7},{11 \over 16}\bigr]
$$
\newline > guess(r);
$$
{n+3\over 2n}
$$
}} 
\medskip\noindent
\proclaim Theorem \sec.7.
We have
$$
\sum_{T\in\setB}x^{|T|} \prod_{h\in \setH(T)} {h+3\over 2h}
= \bigl({1-\sqrt{1-4x}\over 2x}\bigr)^2. 
\leqno{(\sec.13)}
$$
More generally,
$$
\sum_{T\in\setB}x^{|T|} \prod_{h\in \setH(T)} {
\prod_{i=1}^{h-1} (z+2h-i) \over
2h \prod_{i=1}^{h-2} (2z+2h-2-i)     }
= \bigl({1-\sqrt{1-4x}\over 2x}\bigr)^z,
\leqno{(\sec.14)}
$$
or
$$
\sum_{T\in\setB(n)} \prod_{h\in \setH(T)} {
\prod_{i=1}^{h-1} (z+2h-i) \over
2h \prod_{i=1}^{h-2} (2z+2h-2-i)     }
= {z\over n!}\prod_{i=1}^{n-1}(2n-i+z). \leqno{(\sec.15)}
$$

Theorem \sec.5, \sec.6 and \sec.7 
can be derived from the next Theorem by taking
$a=1$, $a=0$ and
$a\rightarrow \infty$ respectively.

\runhookexp[4pt]{18mm}{\hbox{\tt (M.\sec.8)}}{{
\newline > [seq(a+1/n, n=1..7)];  
$$
\bigl[a+1,\; a+{1\over 2},\;a+{1\over 3},\;a+{1\over 4},\;
a+{1\over 5},\;a+{1\over 6},\;a+{1\over 7}\bigr]
$$
\newline > f:=hookgen(\%): hookexp(f\^{}z, 5): map(factor, \%);
$$
\leqalignno{
& \bigl[z(a+1),\; {za+3a+z+1 \over 4},\; 
{(za+5a+z+1)(za+4a+z+2) \over 18a+12za+6+12z},\cr
& {(za+z+2+6a)(za+z+1+7a)(za+z+3+5a) \over 16(2za+5a+1+2z)(za+2a+z+1)},\cr
& {(za+z+4+6a)(za+z+1+9a)(za+z+3+7a)(za+z+2+8a) \over 
20(2za+2z+3+5a)(za+3a+z+1)(2za+2z+1+7a)} 
\bigr]\cr
}
$$
}} 
\medskip\noindent
\proclaim Theorem \sec.8.
Let
$$
\leqalignno{
\sum_{T\in\setB(n)} & \prod_{h\in \setH(T)} {
\prod_{i=1}^{h-1} (za+z+(2h-i)a+i) \over
2h \prod_{i=1}^{h-2} (2za+2z+(2h-2-i)a+i)     }\cr
&= {z(a+1)\over n!}\prod_{i=1}^{n-1}(za+z+(2n-i)a+i).&(\sec.16)\cr
}
$$

{\it Proof}.
Let $f(x)$ be a power series in $x$ defined by
$$
f(x)=1+(a -1) x f(x)^{2a/(a-1)}
$$
and $U_n(z,a)$ be the right-hand side of (\sec.16).
Let $g(x)=f(x)-1$, then
$$
g(x)=(a -1) x (g(x)+1)^{2a/(a-1)}.
$$
By the Lagrange inversion formula we have 
$$
\leqalignno{
[x^n] (g(x)+1)^{z} 
&=
{1\over n} [x^{n-1}] \Bigl( 
z(x+1)^{z-1} (a-1)^n (x+1)^{2an/(a-1)}
\Bigr) \cr
&= {z(a-1)^n\over n} [x^{n-1}] (x+1)^{z-1+2an/(a-1)} \cr
&= {z(a-1)^n\over n (n-1)!} \prod_{i=0}^{n-2}\bigl(z-1+{2an\over a-1} -i\bigr). \cr
&= {z(a-1)\over n!} \prod_{i=1}^{n-1}\bigl(z(a-1)+{2an } -i(a-1)\bigr), \cr
}
$$
so that
$$
[x^n] f(x)^{z((a+1)/(a-1)} 
= {z(a+1)\over n!} \prod_{i=1}^{n-1}\bigl(z(a+1)+{2an } -i(a-1)\bigr)
=U_n(z,a).
$$
By (\sec.3) we have
$$
\rho(n)={[x^n]  f(x)^{z (a+1)/(a-1)}  \over [x^{n-1}] f(x)^{2z (a+1)/(a-1)} }
={U_n(z,a) \over U_{n-1}(2z,a)}.\qed
$$

{\it Remark}. 
Theorem \sec.8 unifies a lot of formulas, including 
the two classical hook formulas (6.4) and (6.5), Postnikov's formula (6.7)
and the generalization due to Lascoux, Du and Liu, another generalization
of Postnikov's formula (6.10), Theorems 6.6 and 6.7.

\runhookexp[4pt]{18mm}{\hbox{\tt (M.\sec.9)}}{{
\newline > hookexp(tan(x)+sec(x), 8);
$$
\bigl[1,{1 \over 4},{1 \over 6},{1 \over 8},{1 \over 10},
{1 \over 12},{1 \over 14},{1 \over 16}\bigr]
$$
\newline > hookexp(z*tan(x)+sec(x), 8);
$$
\bigl[z,{1 \over 4z},{z \over 3+3z^2},{1 \over 8z},
{z \over 5+5z^2},{1 \over 12z},{z \over 7+7z^2},{1 \over 16z}\bigr]
$$
}} 
\medskip\noindent
\proclaim Theorem \sec.9.
We have
$$
\sum_{T\in\setB}x^{|T|} \ \prod_{h\in \setH(T), h\geq 2} \ 
{1\over 2h}
= \tan(x)+\secan(x)
\leqno{(\sec.17)}
$$
and
$$
\sum_{T\in\setB}x^{|T|} \ \prod_{h\in \setH(T)} \rho(h) 
= z\tan(x)+\secan(x),
\leqno{(\sec.18)}
$$
where 
$$
\rho(n)=\cases{
z, &if $n=1$;\cr
\noalign{\smallskip}
{1\over 2nz}, &if $n$ is even; \cr
\noalign{\smallskip}
{z\over n(1+z^2)}, &if $n\geq 3$ is odd.\cr
}
$$

{\it Proof}.
Let $f(x)=z\tan(x)+\secan(x)$. Then
$$
f^2(x)= {1+z^2\sin^2(x)\over \cos^2(x) } + {2z\sin(x)\over \cos^2(x)}.
$$

It is easy to verify that $\rho(1)=z$. For each $k\geq 1$,
$$
\leqalignno{
[x^{2k}]\ f^2(x) 
& = [x^{2k}]\ {1+z^2 \sin^2(x)\over \cos^2(x)} -1 
= [x^{2k}]\ {1+z^2 \sin^2(x)-\cos^2(x)\over \cos^2(x)}  \cr
&= (1+z^2) [x^{2k}]\ \tan^2(x) 
= {1+z^2\over z} [x^{2k}]\ z\tan^2(x) \cr
&= {1+z^2\over z} [x^{2k}]\ (z\tan(x))' -z) \cr 
&= {1+z^2\over z} (2k+1)[x^{2k+1}]\ (z\tan(x)).  \cr
}
$$
So that
$$
\rho(2k+1) = {[x^{2k+1}] f(x) \over [x^{2k}] f^2(x)} ={z\over (2k+1) (1+z^2)}.
$$
On the other hand,
$$
[x^{2k-1}] {2z\sin(x)\over \cos^2(x)}
=2z[x^{2k-1}] {\secan(x)'}
=2z(2k) [x^{2k}] {\secan(x)},
$$
hence
$$
\rho(2k) = {[x^{2k}] f(x) \over [x^{2k-1}] f^2(x)} ={1\over 2z(2k)}.\qed
$$

{\it Remark}. Formula (\sec.17) has a combinatorial interpretation 
due to Foata, Sch\"u\-tzen\-berger and Strehl [FS73; FS74; Vi81; FH01] 
by using the model of  Andr\'e permutations. Note the difference with Theorem 7.1.

\runhookexp[4pt]{18mm}{\hbox{\tt (M.\sec.10)}}{{
\newline > hookexp((1+x)/(1+x\^{}2), 9);
$$
\bigl[1,\;{-1 \over 2},\;1,\;{-1 \over 4},\;1,\;{-1 \over 6},\;1,\;
{-1 \over 8},\;1\bigr]
$$
}} 
\medskip\noindent
\proclaim Theorem \sec.10.
We have
$$
\sum_{T\in\setB}x^{|T|} \ \prod_{h\in \setH(T), h {\rm\ even}} {-1\over h} 
= {1+x\over 1+x^2}.
\leqno{(\sec.19)}
$$

{\it Proof.}
Let
$$
f(x)={1\over 1+x^2} + {x\over 1+x^2}.
$$
Then 
$$
f(x)^2 ={1\over 1+x^2} + {2x\over (1+x^2)^2}.
$$
We have
$$
[x^{2k-1}] {2x\over (1+x^2)^2} 
=-[x^{2k-1}] \bigl( {1\over 1+x^2} \bigr)'
=-(2k)[x^{2k}] {1\over 1+x^2},
$$
so that
$$
\rho(2k) = {[x^{2k}] f(x) \over [x^{2k-1}] f^2(x)} =-{1\over 2k}
$$
and
$$
\rho(2k+1) = {[x^{2k+1}] f(x) \over [x^{2k}] f^2(x)} =1.\qed
$$

\runhookexp[4pt]{18mm}{\hbox{\tt (M.\sec.11)}}{{
\newline > hookexp((1+x)/(1+x\^{}3), 12);
$$
\bigl[1,\;0,\;-1,\;{1 \over 2},\;0,\;{-1 \over 2},\;{1 \over 3},\;0,\;
{-1 \over 3},\; {1 \over 4},\;0,\;{-1 \over 4}\bigr]
$$
}} 
\medskip\noindent
\proclaim Theorem \sec.11.
We have
$$
\sum_{T\in\setB}x^{|T|} \ \prod_{h\in \setH(T)} \rho(h) 
= {1+x\over 1+x^3} 
\leqno{(\sec.20)}
$$
where 
$$
\rho(n)=\cases{
{1/k}; &if $n=3k-2$,\cr
{0}; &if $n=3k-1$,\cr
-{1/k}; &if $n=3k$. \cr
}
$$

{\it Proof.}
Let
$$
f(x)={1\over 1+x^3} + {x\over 1+x^3}.
$$
Then
$$
f(x)^2 ={1\over (1+x^3)^2} + {2x\over (1+x^3)^2}+{x^2\over (1+x^3)^2}.
$$
Thus
$$
\leqalignno{
[x^{3k-1}] {x^2\over (1+x^3)^2} 
&=-{1\over 3}[x^{3k-1}] \bigl( {x\over 1+x^3} \bigr)'
=-k[x^{3k}]  {1\over 1+x^3}. &(\sec.21)\cr
}
$$
It is easy to see that $\rho(3k-1)=0$ and
$$
\rho(3k) = {[x^{3k}] f(x) \over [x^{3k-1}] f^2(x)} =-{1\over k}.
$$
On the other hand, 
$$
[x^{3k-2}] {x\over 1+x^3} = 
[x^{3k-3}] {1\over 1+x^3} =
-[x^{3k}] {1\over 1+x^3}
$$
and by (\sec.25)
$$
[x^{3k-3}] {1\over (1+x^3)^2} = [x^{3k-1}] {x^2\over (1+x^3)^2}
=-k[x^{3k}]  {1\over 1+x^3}. 
$$
Finally,
$$
\rho(3k-2) 
= {[x^{3k-2}] f(x) \over [x^{3k-3}] f^2(x)} ={1\over k}.\qed
$$

Consider the weight function $\rho$ that counts the leaves 
of binary trees. 
\runhookexp[4pt]{18mm}{\hbox{\tt (M.\sec.12)}}{{ 
\newline > [1, seq(2, i=1..7)]; 
$$
[1, 2, 2, 2, 2, 2, 2, 2]
$$
\newline > hookgen(\%);
$$
1+x+4x^2+18x^3+88x^4+456x^5+2464x^6+13736x^7+78432x^8
$$
}} 
\medskip\noindent
The above generating function corresponds to the
sequence A068764 in the on-line encyclopedia of integer sequences [Slo].
It is equal to the generating function of the 
{\it generalized Catalan numbers}. 
\proclaim Theorem \sec.12. 
We have 
$$
\sum_{\l\in\setB} x^{|\l|} \prod_{h\in \setH(\l), h\geq 2}  2
= 
{1-\sqrt{1-8x(1-x)}\over 4x}.
\leqno{(\sec.22)}
$$

More generally,
\runhookexp[4pt]{18mm}{\hbox{\tt (M.\sec.13)}}{{ 
\newline > [z, seq(1, i=1..7)]; 
$$
[z, 1, 1, 1, 1, 1, 1, 1]
$$
\newline > hookgen(\%);
$$
\leqalignno{
&1+(z)x+(2z)x^2+(2z+z^2)x^3+(3z+2z^2)x^4+(2z+5z^2)x^5\cr
&+(4z+6z^2+z^3)x^6+(2z+11z^2+2z^3)x^7+(4z+13z^2+5z^3)x^8\cr
}
$$
}} 
\medskip\noindent
The coefficient of $x^n z^j$ in the above generating function is the 
number of binary trees with $n$ vertices and $j$ leaves, it is given by
the following formulas [Pr96].
\proclaim Theorem \sec.13 [Prodinger]. 
We have 
$$
\sum_{\l\in\setB} x^{|\l|} \prod_{h\in \setH(\l), h =1}  z
= 
\sum_{n\geq 0,j\geq 1}A_{n,j}x^nz^j,
\leqno{(\sec.17)}
$$
where
$$
A_{n,j} = {
2^{n+1-2j} (n-1)!
\over 
j!(j-1)!(n+1-2j)!
}.
$$

\def\sec{7}
\section{\sec. Hook length formulas for complete binary trees} 

A {\it complete binary tree} $T$ is a binary tree such that the two
subtrees of each vertex $v$ are either both empty or both non-empty,
except when $v$ is the latest vertex in the so-called 
{\it inorder} [Kn98a,p.319; Vi81].
For example, there are five complete binary trees with $n=7$ vertices.
%

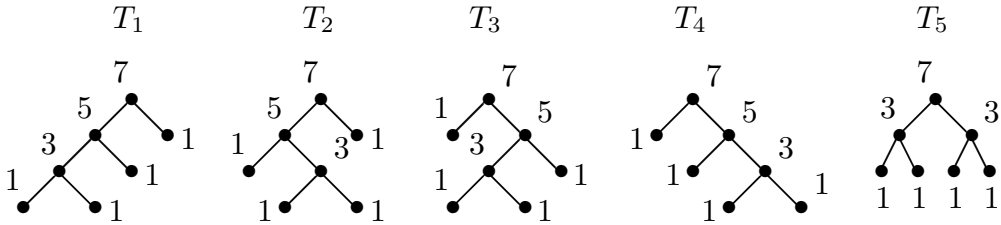
\long\def\maplebegin#1\mapleend{}

\maplebegin
# --------------- begin maple ----------------------

# Copy the following text  to "makefig.mpl"
# then in maple > read("makefig.mpl");
# it will create a file "z_fig_by_maple.tex"

#\unitlength=1pt

Hu:= 4; # height quantities

X0:=5.0; Y0:=0; # origin position

File:=fopen("z_fig_by_maple.tex", WRITE);

pp:=proc(x,y) # X0,Y0 = offset
local a,b,len;
	len:=1;
	fprintf(File, "\\pline(
end;

nn:=proc(x,y) # X0,Y0 = offset
local a,b,len;
	len:=1;
	fprintf(File, "\\nline(
end;

ppe:=proc(x,y) # X0,Y0 = offset
local a,b,len;
	len:=1;
	fprintf(File, "\\plinee(
end;

nne:=proc(x,y) # X0,Y0 = offset
local a,b,len;
	len:=1;
	fprintf(File, "\\nlinee(
end;

mydot:=proc(x,y) # X0,Y0 = offset
local a,b,len;
	len:=1;
	fprintf(File, "\\mydot(
end;

mylabel:=proc(x,y, text) # X0,Y0 = offset
local a,b,len;
	len:=1;
	fprintf(File, "\\mylabel(
end;

mylabel2:=proc(x,y, text) # X0,Y0 = offset
local a,b,len;
	len:=1;
	fprintf(File, "\\mylabel(
end;

mylabel3:=proc(x,y, text) # X0,Y0 = offset
local a,b,len;
	len:=1;
	fprintf(File, "\\mylabel(
end;

mylabel4:=proc(x,y, text) # X0,Y0 = offset
local a,b,len;
	len:=1;
	fprintf(File, "\\mylabel(
end;

mylabel5:=proc(x,y, text) # X0,Y0 = offset
local a,b,len;
	len:=1;
	fprintf(File, "\\mylabel(
end;

lab:=proc(x,y,t) mydot(x,y); mylabel(x,y,t); end;
lab2:=proc(x,y,t) mydot(x,y); mylabel2(x,y,t); end;
lab3:=proc(x,y,t) mydot(x,y); mylabel3(x,y,t); end;
lab4:=proc(x,y,t) mydot(x,y); mylabel4(x,y,t); end;
lab5:=proc(x,y,t) mydot(x,y); mylabel5(x,y,t); end;

DXX:=20;

X0:=-2.8*DXX;
pp(1,0); pp(2,1); pp(3,2); nn(2,1); nn(3,2); nn(4,3);
lab(2,1, "3"); lab(3,2, "5"); lab(4,3, "7");
lab(1,0, "1"); lab3(3,0, "1"); lab3(4,1, "1"); 
lab3(5,2, "1");
mylabel(4,4.5, "$T_1$");

X0:=X0+DXX+DXX/4;

nn(2,2); pp(2,2); pp(1,1); nn(3,3); pp(2,0); nn(3,1);
lab2(3,1, "3"); lab(2,2, "5"); lab(3,3, "7");
lab(1,1,"1"); lab4(2,0,"1"); lab3(4,0,"1"); lab3(4,2,"1");
mylabel(3,4.5, "$T_2$");

X0:=X0+DXX+DXX/8;
nn(2,3); pp(2,1); pp(1,2); pp(1,0); nn(3,2); nn(2,1); 
lab(2,1, "3"); lab2(3,2, "5"); lab2(2,3, "7");
lab(1,2, "1"); lab(1,0, "1"); lab3(3,0, "1"); lab3(4,1,"1");
mylabel(2,4.5, "$T_3$");

X0:=X0+DXX+DXX/3; 

nn(1,3); nn(2,2); pp(0,2); pp(1,1); pp(2,0); nn(3,1);
lab2(3,1, "3"); lab2(2,2, "5"); lab2(1,3, "7");
lab4(0,2,"1"); lab4(1,1,"1"); lab4(2,0,"1"); lab2(4,0,"1");
mylabel(1,4.5, "$T_4$");

X0:=X0+DXX+DXX/7;
pp(1,2); nn(2,3); ppe(0+0.5,1); ppe(2+0.5,1); nne(1,2); nne(3,2);
lab(1,2, "3"); lab2(3,2, "3"); lab(2,3, "7");
lab5(0.5,1, "1");
lab5(1.5,1, "1");
lab5(2.5,1, "1");
lab5(3.5,1, "1");
mylabel(2,4.5, "$T_5$");

fclose(File);

# -------------------- end maple -------------------------

\mapleend


\newbox\boxarbre
\def\pline(#1,#2)#3|{\leftput(#1,#2){\lline(1,1){#3}}}
\def\plinee(#1,#2)#3|{\leftput(#1,#2){\lline(1,2){#3}}}
\def\nline(#1,#2)#3|{\leftput(#1,#2){\lline(1,-1){#3}}}
\def\nlinee(#1,#2)#3|{\leftput(#1,#2){\lline(1,-2){#3}}}
\def\mydot(#1,#2)|{\leftput(#1,#2){$\bullet$}}
\def\mylabel(#1,#2)#3|{\leftput(#1,#2){#3}}
\setbox\boxarbre=\vbox{\vskip
24mm\offinterlineskip 
%
\pline(-52.0,0.0)4.0|
\pline(-48.0,4.0)4.0|
\pline(-44.0,8.0)4.0|
\nline(-48.0,4.0)4.0|
\nline(-44.0,8.0)4.0|
\nline(-40.0,12.0)4.0|
\mydot(-48.8,3.2)|
\mylabel(-50.0,6.0){3}|
\mydot(-44.8,7.2)|
\mylabel(-46.0,10.0){5}|
\mydot(-40.8,11.2)|
\mylabel(-42.0,14.0){7}|
\mydot(-52.8,-0.8)|
\mylabel(-54.0,2.0){1}|
\mydot(-44.8,-0.8)|
\mylabel(-42.5,-1.5){1}|
\mydot(-40.8,3.2)|
\mylabel(-38.5,2.5){1}|
\mydot(-36.8,7.2)|
\mylabel(-34.5,6.5){1}|
\mylabel(-42.0,20.0){$T_1$}|
\nline(-23.0,8.0)4.0|
\pline(-23.0,8.0)4.0|
\pline(-27.0,4.0)4.0|
\nline(-19.0,12.0)4.0|
\pline(-23.0,0.0)4.0|
\nline(-19.0,4.0)4.0|
\mydot(-19.8,3.2)|
\mylabel(-17.5,5.5){3}|
\mydot(-23.8,7.2)|
\mylabel(-25.0,10.0){5}|
\mydot(-19.8,11.2)|
\mylabel(-21.0,14.0){7}|
\mydot(-27.8,3.2)|
\mylabel(-29.0,6.0){1}|
\mydot(-23.8,-0.8)|
\mylabel(-26.0,-1.5){1}|
\mydot(-15.8,-0.8)|
\mylabel(-13.5,-1.5){1}|
\mydot(-15.8,7.2)|
\mylabel(-13.5,6.5){1}|
\mylabel(-21.0,20.0){$T_2$}|
\nline(-0.5,12.0)4.0|
\pline(-0.5,4.0)4.0|
\pline(-4.5,8.0)4.0|
\pline(-4.5,0.0)4.0|
\nline(3.5,8.0)4.0|
\nline(-0.5,4.0)4.0|
\mydot(-1.2,3.2)|
\mylabel(-2.5,6.0){3}|
\mydot(2.8,7.2)|
\mylabel(5.0,9.5){5}|
\mydot(-1.2,11.2)|
\mylabel(1.0,13.5){7}|
\mydot(-5.2,7.2)|
\mylabel(-6.5,10.0){1}|
\mydot(-5.2,-0.8)|
\mylabel(-6.5,2.0){1}|
\mydot(2.8,-0.8)|
\mylabel(5.0,-1.5){1}|
\mydot(6.8,3.2)|
\mylabel(9.0,2.5){1}|
\mylabel(-2.5,20.0){$T_3$}|
\nline(22.2,12.0)4.0|
\nline(26.2,8.0)4.0|
\pline(18.2,8.0)4.0|
\pline(22.2,4.0)4.0|
\pline(26.2,0.0)4.0|
\nline(30.2,4.0)4.0|
\mydot(29.4,3.2)|
\mylabel(31.7,5.5){3}|
\mydot(25.4,7.2)|
\mylabel(27.7,9.5){5}|
\mydot(21.4,11.2)|
\mylabel(23.7,13.5){7}|
\mydot(17.4,7.2)|
\mylabel(15.2,6.5){1}|
\mydot(21.4,3.2)|
\mylabel(19.2,2.5){1}|
\mydot(25.4,-0.8)|
\mylabel(23.2,-1.5){1}|
\mydot(33.4,-0.8)|
\mylabel(35.7,1.5){1}|
\mylabel(20.2,20.0){$T_4$}|
\pline(45.0,8.0)4.0|
\nline(49.0,12.0)4.0|
\plinee(43.0,4.0)2.0|
\plinee(51.0,4.0)2.0|
\nlinee(45.0,8.0)2.0|
\nlinee(53.0,8.0)2.0|
\mydot(44.3,7.2)|
\mylabel(43.0,10.0){3}|
\mydot(52.3,7.2)|
\mylabel(54.5,9.5){3}|
\mydot(48.3,11.2)|
\mylabel(47.0,14.0){7}|
\mydot(42.3,3.2)|
\mylabel(42.5,0.0){1}|
\mydot(46.3,3.2)|
\mylabel(46.5,0.0){1}|
\mydot(50.3,3.2)|
\mylabel(50.5,0.0){1}|
\mydot(54.3,3.2)|
\mylabel(54.5,0.0){1}|
\mylabel(47.0,20.0){$T_5$}|%
}
$$
\kern-4mm\box\boxarbre
$$
\medskip\noindent
We have the hook length multi-sets 
$\setH(T_1)=\setH(T_2)= \setH(T_3)= \setH(T_4)= \{1,1,1,1,3,5,7\}$
and $\setH(T_5)=\{1,1,1,1,3,3,7\}$. 
There are also five complete binary trees with $n=6$ vertices.
%

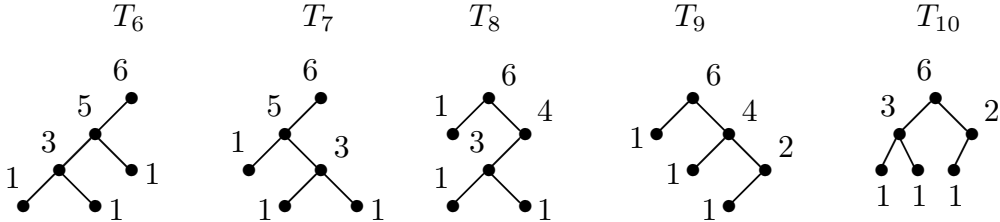
\long\def\maplebegin#1\mapleend{}

\maplebegin

# --------------- begin maple ----------------------

# Copy the following text  to "makefig.mpl"
# then in maple > read("makefig.mpl");
# it will create a file "z_fig_by_maple.tex"

#\unitlength=1pt

Hu:= 4; # height quantities

X0:=5.0; Y0:=0; # origin position

File:=fopen("z_fig_by_maple.tex", WRITE);

pp:=proc(x,y) # X0,Y0 = offset
local a,b,len;
	len:=1;
	fprintf(File, "\\pline(
end;

nn:=proc(x,y) # X0,Y0 = offset
local a,b,len;
	len:=1;
	fprintf(File, "\\nline(
end;

ppe:=proc(x,y) # X0,Y0 = offset
local a,b,len;
	len:=1;
	fprintf(File, "\\plinee(
end;

nne:=proc(x,y) # X0,Y0 = offset
local a,b,len;
	len:=1;
	fprintf(File, "\\nlinee(
end;

mydot:=proc(x,y) # X0,Y0 = offset
local a,b,len;
	len:=1;
	fprintf(File, "\\mydot(
end;

mylabel:=proc(x,y, text) # X0,Y0 = offset
local a,b,len;
	len:=1;
	fprintf(File, "\\mylabel(
end;

mylabel2:=proc(x,y, text) # X0,Y0 = offset
local a,b,len;
	len:=1;
	fprintf(File, "\\mylabel(
end;

mylabel3:=proc(x,y, text) # X0,Y0 = offset
local a,b,len;
	len:=1;
	fprintf(File, "\\mylabel(
end;

mylabel4:=proc(x,y, text) # X0,Y0 = offset
local a,b,len;
	len:=1;
	fprintf(File, "\\mylabel(
end;

mylabel5:=proc(x,y, text) # X0,Y0 = offset
local a,b,len;
	len:=1;
	fprintf(File, "\\mylabel(
end;

lab:=proc(x,y,t) mydot(x,y); mylabel(x,y,t); end;
lab2:=proc(x,y,t) mydot(x,y); mylabel2(x,y,t); end;
lab3:=proc(x,y,t) mydot(x,y); mylabel3(x,y,t); end;
lab4:=proc(x,y,t) mydot(x,y); mylabel4(x,y,t); end;
lab5:=proc(x,y,t) mydot(x,y); mylabel5(x,y,t); end;

DXX:=20;

X0:=-2.8*DXX;
pp(1,0); pp(2,1); pp(3,2); nn(2,1); nn(3,2);
lab(2,1, "3"); lab(3,2, "5"); lab(4,3, "6");
lab(1,0, "1"); lab3(3,0, "1"); lab3(4,1, "1"); 
mylabel(4,4.5, "$T_6$");

X0:=X0+DXX+DXX/4;

nn(2,2); pp(2,2); pp(1,1);  pp(2,0); nn(3,1);
lab2(3,1, "3"); lab(2,2, "5"); lab(3,3, "6");
lab(1,1,"1"); lab4(2,0,"1"); lab3(4,0,"1"); 
mylabel(3,4.5, "$T_7$");

X0:=X0+DXX+DXX/8;
nn(2,3); pp(2,1); pp(1,2); pp(1,0); nn(2,1); 
lab(2,1, "3"); lab2(3,2, "4"); lab2(2,3, "6");
lab(1,2, "1"); lab(1,0, "1"); lab3(3,0, "1"); 
mylabel(2,4.5, "$T_8$");

X0:=X0+DXX+DXX/3; 

nn(1,3); nn(2,2); pp(0,2); pp(1,1); pp(2,0); 
lab2(3,1, "2"); lab2(2,2, "4"); lab2(1,3, "6");
lab4(0,2,"1"); lab4(1,1,"1"); lab4(2,0,"1"); 
mylabel(1,4.5, "$T_9$");

X0:=X0+DXX+DXX/7;
pp(1,2); nn(2,3); ppe(0+0.5,1); ppe(2+0.5,1); nne(1,2); 
lab(1,2, "3"); lab2(3,2, "2"); lab(2,3, "6");
lab5(0.5,1, "1");
lab5(1.5,1, "1");
lab5(2.5,1, "1");
mylabel(2,4.5, "$T_{10}$");

fclose(File);

# -------------------- end maple -------------------------

\mapleend


\newbox\boxarbre
\def\pline(#1,#2)#3|{\leftput(#1,#2){\lline(1,1){#3}}}
\def\plinee(#1,#2)#3|{\leftput(#1,#2){\lline(1,2){#3}}}
\def\nline(#1,#2)#3|{\leftput(#1,#2){\lline(1,-1){#3}}}
\def\nlinee(#1,#2)#3|{\leftput(#1,#2){\lline(1,-2){#3}}}
\def\mydot(#1,#2)|{\leftput(#1,#2){$\bullet$}}
\def\mylabel(#1,#2)#3|{\leftput(#1,#2){#3}}
\setbox\boxarbre=\vbox{\vskip
24mm\offinterlineskip 
%
\pline(-52.0,0.0)4.0|
\pline(-48.0,4.0)4.0|
\pline(-44.0,8.0)4.0|
\nline(-48.0,4.0)4.0|
\nline(-44.0,8.0)4.0|
\mydot(-48.8,3.2)|
\mylabel(-50.0,6.0){3}|
\mydot(-44.8,7.2)|
\mylabel(-46.0,10.0){5}|
\mydot(-40.8,11.2)|
\mylabel(-42.0,14.0){6}|
\mydot(-52.8,-0.8)|
\mylabel(-54.0,2.0){1}|
\mydot(-44.8,-0.8)|
\mylabel(-42.5,-1.5){1}|
\mydot(-40.8,3.2)|
\mylabel(-38.5,2.5){1}|
\mylabel(-42.0,20.0){$T_6$}|
\nline(-23.0,8.0)4.0|
\pline(-23.0,8.0)4.0|
\pline(-27.0,4.0)4.0|
\pline(-23.0,0.0)4.0|
\nline(-19.0,4.0)4.0|
\mydot(-19.8,3.2)|
\mylabel(-17.5,5.5){3}|
\mydot(-23.8,7.2)|
\mylabel(-25.0,10.0){5}|
\mydot(-19.8,11.2)|
\mylabel(-21.0,14.0){6}|
\mydot(-27.8,3.2)|
\mylabel(-29.0,6.0){1}|
\mydot(-23.8,-0.8)|
\mylabel(-26.0,-1.5){1}|
\mydot(-15.8,-0.8)|
\mylabel(-13.5,-1.5){1}|
\mylabel(-21.0,20.0){$T_7$}|
\nline(-0.5,12.0)4.0|
\pline(-0.5,4.0)4.0|
\pline(-4.5,8.0)4.0|
\pline(-4.5,0.0)4.0|
\nline(-0.5,4.0)4.0|
\mydot(-1.2,3.2)|
\mylabel(-2.5,6.0){3}|
\mydot(2.8,7.2)|
\mylabel(5.0,9.5){4}|
\mydot(-1.2,11.2)|
\mylabel(1.0,13.5){6}|
\mydot(-5.2,7.2)|
\mylabel(-6.5,10.0){1}|
\mydot(-5.2,-0.8)|
\mylabel(-6.5,2.0){1}|
\mydot(2.8,-0.8)|
\mylabel(5.0,-1.5){1}|
\mylabel(-2.5,20.0){$T_8$}|
\nline(22.2,12.0)4.0|
\nline(26.2,8.0)4.0|
\pline(18.2,8.0)4.0|
\pline(22.2,4.0)4.0|
\pline(26.2,0.0)4.0|
\mydot(29.4,3.2)|
\mylabel(31.7,5.5){2}|
\mydot(25.4,7.2)|
\mylabel(27.7,9.5){4}|
\mydot(21.4,11.2)|
\mylabel(23.7,13.5){6}|
\mydot(17.4,7.2)|
\mylabel(15.2,6.5){1}|
\mydot(21.4,3.2)|
\mylabel(19.2,2.5){1}|
\mydot(25.4,-0.8)|
\mylabel(23.2,-1.5){1}|
\mylabel(20.2,20.0){$T_9$}|
\pline(45.0,8.0)4.0|
\nline(49.0,12.0)4.0|
\plinee(43.0,4.0)2.0|
\plinee(51.0,4.0)2.0|
\nlinee(45.0,8.0)2.0|
\mydot(44.3,7.2)|
\mylabel(43.0,10.0){3}|
\mydot(52.3,7.2)|
\mylabel(54.5,9.5){2}|
\mydot(48.3,11.2)|
\mylabel(47.0,14.0){6}|
\mydot(42.3,3.2)|
\mylabel(42.5,0.0){1}|
\mydot(46.3,3.2)|
\mylabel(46.5,0.0){1}|
\mydot(50.3,3.2)|
\mylabel(50.5,0.0){1}|
\mylabel(47.0,20.0){$T_{10}$}|
}
$$
\kern-4mm\box\boxarbre
$$
\medskip\noindent
We have 
$\setH(T_6)=\setH(T_7)= \{1,1,1,3,5,6\}$,
$\setH(T_8)= \{1,1,1,3,4,6\}$,
$\setH(T_9)= \{1,1,1,2,4,6\}$
and $\setH(T_{10})=\{1,1,1,2,3,6\}$. 

\medskip
Let $\setC$ (resp. $\setC(n)$) denote
the set of all complete binary trees (resp. all complete 
binary trees with $n$ vertices),
so that 
$$\setC=\bigcup_{n\geq 0} \setC(n).$$
Again, define the {\it hook
length expansion} for complete binary trees by
$$
\sum_{T\in\setC} x^{|T|} \prod_{h\in \setH(T)} \rho(h) 
= f(x),  \leqno{(\sec.1)}
$$
where $f(x)\in K[[x]]$ is a power series in $x$ with coefficients in $K$ such 
that $f(0)=1$. See Sections 5 and 6 for more comments about the hook length
expansion.
Let $f(x)=1+f_1x +f_2x^2+f_3x^3+\cdots$ be the generating function for 
complete binary trees by the weight function  $\rho$.
With each $T\in\setC(n)$ ($n\geq 1$) we can 
associate a triplet $(T', T'', v)$, where
$T'\in\setC(k)$ ($0\leq k\leq n-1$ and $k$ is an odd integer),
$T''\in\setC(n-1-k)$  and the root $v$ of $T$ whose hook length $h_v=n$.
Hence, (\sec.1) is equivalent to
$$
\rho(n)\sum_{k=0, k {\rm\ odd}}^{n-1} f_k f_{n-1-k}=f_n \quad (n\geq 1).
\leqno{(\sec.2)}
$$
Formula (\sec.2) can be used to calculate $f(x)$ for a given $\rho$,
or to calculate $\rho$ for a given $f(x)$.  
It also has the equivalent form
$$
\rho(n)= {[x^n] f(x) \over [x^{n-1}] (f(x)-f(-x))f(x)/2},
\leqno{(\sec.3)}
$$
because
$$
\sum_{k\geq 1, k {\rm\ odd}} f_k x^k = (f(x)-f(-x))/2.
$$

Next we use the maple package {\tt HookExp} to find hook formulas for 
complete binary trees, whose proofs are always based on (\sec.3).
\runhookexp[4pt]{18mm}{\hbox{\tt (M.\sec.1)}}{{
\newline > hooktype:="CBT":  
\newline  \qquad     \# working on complete binary trees
\newline > hookexp(tan(x)+sec(x), 9);  
$$
\bigl[1,{1 \over 2},{1 \over 3},{1 \over 4},{1 \over 5},{1 \over 6},
{1 \over 7},{1 \over 8},{1 \over 9}\bigr]
$$
\newline > hookexp(z*tan(x)+sec(x), 9);  
$$
\bigl[z,{1 \over 2z},{1 \over 3z},{1 \over 4z},{1 \over 5z},{1 \over 6z},{1 \over 7z},{1 \over 8z},{1 \over 9z}\bigr]
$$
}} 
\medskip\noindent
\proclaim Theorem \sec.1.
We have
$$
\sum_{T\in\setC}x^{|T|}  \prod_{h\in \setH(T)} {1\over h}= 
\tan(x)+\secan(x). \leqno{(\sec.4)}
$$
$$
\sum_{T\in\setC}x^{|T|}  \prod_{h\in \setH(T), h=1} z
\prod_{h\in \setH(T), h\geq 2} {1\over zh}= 
z\tan(x)+\secan(x). \leqno{(\sec.5)}
$$

{\it Proof}.  By (\sec.3)
$$
\leqalignno{
\rho(n)
&={ [x^n] z\tan(x)+\secan(x) \over
[x^{n-1}] z\tan(x)(z\tan(x)+\secan(x)) }\cr
&={ [x^n] z\tan(x)+\secan(x) \over
z[x^{n-1}] (z\tan(x)+\secan(x))'-z }\cr
&={ [x^n] \tan(x)+\secan(x) \over
z[x^{n}]n (\tan(x)+\secan(x)) }={1\over zn}.\qed\cr
}
$$

{\it Remark.} Recall that the $n$-th Euler number is the coefficient 
of $x^n/n!$
in the expansion of the series $\tan(x)+\secan(x)$ (see, e.g. [Vi81]). 
It is well-known that
$E_n$ is equal to the number of
alternating permutations of order $n$, which, in turn, is equal to the number
of increasing labeled complete binary trees (See Theorem 6.1). So that
$$
\sum_{T\in\setC(n)} n! \prod_{v\in T} {1\over h_v}= E_n.
$$
This gives a combinatorial proof of Theorem \sec.1.
Note that Theorem 6.9 also involves the Euler numbers, but the combinatorial 
argument is totally different.


\runhookexp[4pt]{18mm}{\hbox{\tt (M.\sec.2)}}{{
\newline > hookexp(exp(x), 9);  
$$
\bigl[1,{1 \over 2},{1 \over 6},{1 \over 16},{1 \over 40},{1 \over 96},
{1 \over 224},{1 \over 512},{1 \over 1152}\bigr]
$$
}} 
\medskip\noindent
\proclaim Theorem \sec.2.
We have
$$
\sum_{T\in\setC}x^{|T|}  \prod_{h\in \setH(T), h\geq 2} {1\over h2^{h-2}}= 
e^x. \leqno{(\sec.6)}
$$

{\it Proof}. From (\sec.3)
$$
\leqalignno{
\rho(n)
&={ [x^n] e^x \over [x^{n-1}] (e^x-e^{-x})e^x/2 } 
={2 [x^n] e^x \over [x^{n-1}] e^{2x}-1 } \cr
&={ 2/n! \over 2^{n-1}/(n-1)! } = {1\over n2^{n-2}}.\qed \cr
}
$$

{\it Remark}. We do not have any combinatorial proof of Theorem \sec.2. 
See also [Ha08b]. It can be viewed as a complete binary tree version
of Theorem 6.3.

\runhookexp[4pt]{18mm}{\hbox{\tt (M.\sec.3)}}{{
\newline > hookexp(1/(1-x), 14);  
$$
\bigl[1,1,1,{1 \over 2},{1 \over 2},{1 \over 3},{1 \over 3},
{1 \over 4},{1 \over 4},{1 \over 5},{1 \over 5},{1 \over 6},
{1 \over 6},{1 \over 7}\bigr]
$$
}} 
\medskip\noindent
\proclaim Theorem \sec.3.
We have
$$
\sum_{T\in\setC}x^{|T|}  \prod_{h\in \setH(T), h\geq 2} \rho(h)= 
{1\over 1-x}, \leqno{(\sec.7)}
$$
where
$$
\rho(n)=\cases{
1, &if $n=1$;\cr
1/k, &if $n=2k+1$ ($k\geq 1$);\cr
1/k, &if $n=2k$ ($k\geq 1$).\cr
} \leqno{(\sec.8)}
$$

{\it Proof}. Since
$$
{x\over (1-x)(1-x^2)} =
{x+x^2\over (1-x^2)^2} = (x+x^2) \sum_{k\geq 0} (k+1)x^{2k},
$$
we have
$$
[x^{2k}] {x\over (1-x)(1-x^2)} =
[x^{2k-1}] {x\over (1-x)(1-x^2)} = k. 
$$
By (\sec.3)
$$
\leqalignno{
\rho(n)
&={ [x^n] 1/(1-x) \over [x^{n-1}] (1/(1-x)-1/(1+x))/(1-x)/2 } \cr
&={1 \over [x^{n-1}] x/(1-x)/(1-x^2) }. \qed\cr
}
$$

{\it Remark}. For each complete binary tree $T$ of $2k$ or $2k+1$ vertices
we obtain, in a bijective manner, 
a binary tree $T'$ of $k$ vertices by deleting all leaves of $T$ [Kn98a, p.399].
This gives a combinatorial proof of Theorem~\sec.3 via Theorem~\sec.1.

\runhookexp[4pt]{18mm}{\hbox{\tt (M.\sec.4)}}{{
\newline > hookexp( (1-sqrt(1-4*x\^{}2))/(2*x\^{}2)*(1+x), 9); 
$$
[1, 1, 1, 1, 1, 1, 1]
$$
\newline > hookexp( (1-sqrt(1-4*x\^{}2))/(2*x\^{}2)*(1+z*x), 9); 
$$
\bigl[z,{1 \over z},{1 \over z},{1 \over z},{1 \over z},
{1 \over z},{1 \over z}\bigr]
$$
}} 
\medskip\noindent
\proclaim Theorem \sec.4.
We have
$$
\sum_{T\in\setC}x^{|T|}  \prod_{h\in \setH(T)} 1= 
{1-\sqrt{1-4x^2}\over 2x^2}(1+x)
\leqno{(\sec.9)}
$$
and
$$
\sum_{T\in\setC}x^{|T|}  \prod_{h\in \setH(T), h=1} z  
\prod_{h\in \setH(T), h\geq 2} {1\over z}         = 
{1-\sqrt{1-4x^2}\over 2x^2}(1+zx)
\leqno{(\sec.10)}
$$

{\it Proof}. 
Let $f(x)$ be the right-hand side of (\sec.10). We can verify 
$$
{ (f(x)-f(-x))f(x)\over 2} \times {x\over z} =f(x) -1-zx.
$$
$$
\rho(n)
={ [x^n] f(x) \over [x^{n-1}] (f(x)-f(-x))f(x)/2 } ={1\over z}.\qed
$$

{\it Remark}. The bijection between binary trees and 
complete binary trees described in Theorem \sec.3 gives a combinatorial proof
of Theorem \sec.4 via Theorem 6.2.

\runhookexp[4pt]{18mm}{\hbox{\tt (M.\sec.5)}}{{
\newline > hookexp( (1+x)/(1+x\^{}2), 11);
$$
\bigl[1,-1,-1,{-1 \over 2},{-1 \over 2},{-1 \over 3},{-1 \over 3},
{-1 \over 4},{-1 \over 4},{-1 \over 5},{-1 \over 5}
\bigr]
$$
}} 
\medskip\noindent
\proclaim Theorem \sec.5.
We have
$$
\sum_{T\in\setC}x^{|T|}  \prod_{h\in \setH(T), h\geq 2} \rho(h)= 
{1+x\over 1+x^2}, \leqno{(\sec.11)}
$$
where
$$
\rho(n)=\cases{
1, &if $n=1$;\cr
-1/k, &if $n=2k+1$ ($k\geq 1$);\cr
-1/k, &if $n=2k$ ($k\geq 1$).\cr
} \leqno{(\sec.12)}
$$

{\it Proof}. 
Let $f(x)$ be the right-hand side of (\sec.11)
$$
f(x):= {1+x\over 1+x^2} 
=(1+x) \sum_{k\geq 0} (-1)^k x^{2k}
$$
and
$$
F(x):={ (f(x)-f(-x))f(x)\over 2} = {x+x^2\over (1+x^2)^2} 
=(x+x^2) \sum_{k\geq 0} (k+1)(-1)^k x^{2k}.
$$
We have
$$
\leqalignno{
[x^{2k+1}] f(x) &= [x^{2k}] f(x) = (-1)^{k};\cr
[x^{2k}] F(x) &= [x^{2k-1}] F(x) = (-1)^{k-1}k. \cr
}
$$
By (\sec.3)
$$
\rho(n)
={ [x^n] f(x) \over [x^{n-1}] F(x) }  = -{1\over k}. \qed
$$

\runhookexp[4pt]{18mm}{\hbox{\tt (M.\sec.6)}}{{
\newline > hookexp( (1+x)/(1+x\^{}4), 4);
\newline {\tt Denominator is zero, no solution for n=4.}
\vskip 4pt
\newline > hookexp( (1+x)/(1+x\^{}3), 16);
$$
\bigl[1,\;0,\;-1,\;1,\;0,\;-1,\;1,\;0,\;{-1 \over 2},\;{1 \over 2},\;0,\;
{-1 \over 2},\; {1 \over 2},\;0,\;{-1 \over 3},\;{1 \over 3}\bigr]
$$
}} 
\medskip\noindent
\proclaim Theorem \sec.6.
We have
$$
\sum_{T\in\setC}x^{|T|}  \prod_{h\in \setH(T), h\geq 2} \rho(h)= 
{1+x\over 1+x^3}, \leqno{(\sec.13)}
$$
where
$$
\rho(n)=\cases{
1, &if $n=1$;\cr
0, &if $n=3k-1$ ($k\geq 1$);\cr
-1/k, &if $n=6k-3$ or $n=6k$ ($k\geq 1$);\cr
1/k, &if $n=6k-2$ or $n=6k+1$ ($k\geq 1$).\cr
} \leqno{(\sec.14)}
$$

{\it Proof}. 
Let $f(x)$ be the right-hand side of (\sec.13)
$$
f(x):= {1+x\over 1+x^3} 
=(1+x) \sum_{k\geq 0} (-1)^k x^{3k}
$$
and
$$
\leqalignno{
F(x)
&:={ (f(x)-f(-x))f(x)\over 2}  \cr
& = {x+x^2-x^3-2x^4-x^5+x^6+x^7\over (1-x^6)^2}  \cr
& = (x+x^2-x^3-2x^4-x^5+x^6+x^7) \sum_{k\geq 0} (k+1) x^{6k}.\cr
}
$$
We have
$$ 
\leqalignno{
[x^{3k+2}] f(x) &=0; \cr 
[x^{3k+1}] f(x) &= [x^{3k}] f(x) = (-1)^{k}; \cr
[x^{6k+2}] F(x) &= [x^{6k+6}] F(x)= k+1; \cr
[x^{6k+3}] F(x) &= [x^{6k+5}] F(x)= -(k+1); \cr 
[x^{6k+4}] F(x) &=  -2(k+1). \cr 
}$$
By (\sec.3)
$$
\rho(n)
={ [x^n] f(x) \over [x^{n-1}] F(x) } . \qed
$$

\def\sec{8}
\section{\sec. Hook length formulas for Fibonacci trees} 
A {\it Fibonacci tree} $T$ is a binary tree such that the right 
subtree of each vertex $v$ is a binary tree with 0 or 1 vertex [St75, SY89].
For example, there are five Fibonacci trees with $n=4$ vertices.
%

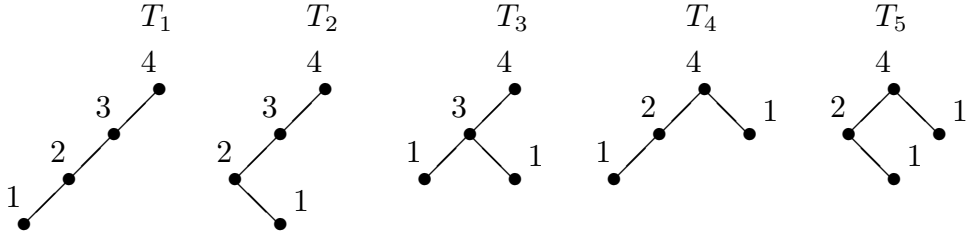
\long\def\maplebegin#1\mapleend{}

\maplebegin

# --------------- begin maple ----------------------

# Copy the following text  to "makefig.mpl"
# then in maple > read("makefig.mpl");
# it will create a file "z_fig_by_maple.tex"

#\unitlength=1pt

Hu:= 5; # height quantities

X0:=0.0; Y0:=0; # origin position

File:=fopen("z_fig_by_maple.tex", WRITE);

pline:=proc(x,y) # X0,Y0 = offset
local a,b,len;
	len:=1;
	fprintf(File, "\\pline(
end;

nline:=proc(x,y) # X0,Y0 = offset
local a,b,len;
	len:=1;
	fprintf(File, "\\nline(
end;

mydot:=proc(x,y) # X0,Y0 = offset
local a,b,len;
	len:=1;
	fprintf(File, "\\mydot(
end;

mylabel:=proc(x,y, text) # X0,Y0 = offset
local a,b,len;
	len:=1;
	fprintf(File, "\\mylabel(
end;

mylabel2:=proc(x,y, text) # X0,Y0 = offset
local a,b,len;
	len:=1;
	fprintf(File, "\\mylabel(
end;

dotlabel:=proc(x,y,t) mydot(x,y); mylabel(x,y,t); end;
dotlabel2:=proc(x,y,t) mydot(x,y); mylabel2(x,y,t); end;

DXX:=21;

X0:=-2.6*DXX;
pline(2,1); pline(3,2); pline(1,0);
dotlabel(2,1, "2"); dotlabel(3,2, "3"); dotlabel(4,3, "4");
dotlabel(1,0,"1");
mylabel(4,4, "$T_1$");

X0:=X0+7/8*DXX;
pline(2,1); pline(3,2); nline(2,1);
dotlabel(2,1, "2"); dotlabel(3,2, "3"); dotlabel(4,3, "4");
dotlabel2(3,0,"1");
mylabel(4,4, "$T_2$");

X0:=X0+DXX;
pline(2,1); pline(3,2); nline(3,2);
dotlabel(2,1, "1"); dotlabel(3,2, "3"); dotlabel(4,3, "4");
dotlabel2(4,1,"1");
mylabel(4,4, "$T_3$");

X0:=X0+DXX;
pline(2,1); pline(3,2); nline(4,3);
dotlabel(2,1, "1"); dotlabel(3,2, "2"); dotlabel(4,3, "4");
dotlabel2(5,2,"1");
mylabel(4,4, "$T_4$");

X0:=X0+DXX;
pline(3,2); nline(4,3); nline(3,2);
dotlabel(3,2, "2"); dotlabel(4,3, "4");
dotlabel2(4,1,"1"); dotlabel2(5,2,"1");
mylabel(4,4, "$T_5$");

fclose(File);

# -------------------- end maple -------------------------

\mapleend


\newbox\boxarbre
\def\pline(#1,#2)#3|{\leftput(#1,#2){\lline(1,1){#3}}}
\def\nline(#1,#2)#3|{\leftput(#1,#2){\lline(1,-1){#3}}}
\def\mydot(#1,#2)|{\leftput(#1,#2){$\bullet$}}
\def\mylabel(#1,#2)#3|{\leftput(#1,#2){#3}}
\setbox\boxarbre=\vbox{\vskip
27mm\offinterlineskip 
%
\pline(-44.6,5.0)5.0|
\pline(-39.6,10.0)5.0|
\pline(-49.6,0.0)5.0|
\mydot(-45.4,4.2)|
\mylabel(-46.6,7.0){2}|
\mydot(-40.4,9.2)|
\mylabel(-41.6,12.0){3}|
\mydot(-35.4,14.2)|
\mylabel(-36.6,17.0){4}|
\mydot(-50.4,-0.8)|
\mylabel(-51.6,2.0){1}|
\mylabel(-36.6,22.0){$T_1$}|
\pline(-26.2,5.0)5.0|
\pline(-21.2,10.0)5.0|
\nline(-26.2,5.0)5.0|
\mydot(-27.0,4.2)|
\mylabel(-28.2,7.0){2}|
\mydot(-22.0,9.2)|
\mylabel(-23.2,12.0){3}|
\mydot(-17.0,14.2)|
\mylabel(-18.2,17.0){4}|
\mydot(-22.0,-0.8)|
\mylabel(-19.7,1.5){1}|
\mylabel(-18.2,22.0){$T_2$}|
\pline(-5.2,5.0)5.0|
\pline(-0.2,10.0)5.0|
\nline(-0.2,10.0)5.0|
\mydot(-6.0,4.2)|
\mylabel(-7.2,7.0){1}|
\mydot(-1.0,9.2)|
\mylabel(-2.2,12.0){3}|
\mydot(4.0,14.2)|
\mylabel(2.8,17.0){4}|
\mydot(4.0,4.2)|
\mylabel(6.3,6.5){1}|
\mylabel(2.8,22.0){$T_3$}|
\pline(15.8,5.0)5.0|
\pline(20.8,10.0)5.0|
\nline(25.8,15.0)5.0|
\mydot(15.0,4.2)|
\mylabel(13.8,7.0){1}|
\mydot(20.0,9.2)|
\mylabel(18.8,12.0){2}|
\mydot(25.0,14.2)|
\mylabel(23.8,17.0){4}|
\mydot(30.0,9.2)|
\mylabel(32.3,11.5){1}|
\mylabel(23.8,22.0){$T_4$}|
\pline(41.8,10.0)5.0|
\nline(46.8,15.0)5.0|
\nline(41.8,10.0)5.0|
\mydot(41.0,9.2)|
\mylabel(39.8,12.0){2}|
\mydot(46.0,14.2)|
\mylabel(44.8,17.0){4}|
\mydot(46.0,4.2)|
\mylabel(48.3,6.5){1}|
\mydot(51.0,9.2)|
\mylabel(53.3,11.5){1}|
\mylabel(44.8,22.0){$T_5$}|
}
$$
\kern-4mm\box\boxarbre
$$
\medskip\noindent
We have the hook length multi-sets 
$\setH(T_1)=\setH(T_2)= \{1,2,3,4\}$,
$\setH(T_4)=\setH(T_5)= \{1,1,2,4\}$
and $\setH(T_3)=\{1,1,3,4\}$. 

\medskip
Let $\setF$ (resp. $\setF(n)$) denote
the set of all Fibonacci trees (resp. all 
Fibonacci trees with $n$ vertices),
so that 
$$\setF=\bigcup_{n\geq 0} \setF(n).$$
As for binary trees, we define the {\it hook length expansion} for 
Fibonacci trees by
$$
\sum_{T\in\setF} x^{|T|} \prod_{h\in \setH(T)} \rho(h) 
= f(x),  \leqno{(\sec.1)}
$$
where $f(x)\in K[[x]]$ is a power series in $x$ with coefficients in $K$ such 
that $f(0)=1$. See Sections 5 and 6 for more comments about the hook length
expansion.
Let $f(x)=1+f_1x +f_2x^2+f_3x^3+\cdots$ be the generating function for  
Fibonacci trees by the weight function  $\rho$.
By definition of Fibonacci trees
formula (\sec.1) is equivalent to
$$
f_n=\rho(n) f_{n-1} + \rho(n)\rho(1) f_{n-2}. \leqno{(\sec.2)} 
$$
Formula (\sec.2) can be used to calculate $f(x)$ for a given $\rho$,
or to calculate $\rho$ for a given $f(x)$.  
Next we use the maple package {\tt HookExp} to find hook formulas for 
Fibonacci trees, whose proofs are always based on (\sec.2).

\runhookexp[4pt]{18mm}{\hbox{\tt (M.\sec.1)}}{{
\newline > hooktype:="FT"  \# working on Fibonacci trees
\newline > hookexp(1/(1-x-x\^{}2), 9); 
$$
[1, 1, 1, 1, 1, 1, 1, 1, 1]
$$
}} 
\smallskip\noindent
\proclaim Theorem \sec.1.
We have
$$
\sum_{T\in\setF}x^{|T|}  \prod_{h\in \setH(T)} {1}= 
{1\over 1-x-x^2}. \leqno{(\sec.3)}
$$

{\it Proof}. Let $f(x)=1+\sum_{n\geq 1}f_n x^n$ be the right-hand side of
(\sec.3), then $f_n=f_{n-1}+f_{n-2}$. Relation (\sec.2) is verified.\qed

{\it Remark}. The number of Fibonacci trees with $n$ vertices is the 
$n$-th Fibonacci number. 
\medskip

\runhookexp[4pt]{18mm}{\hbox{\tt (M.\sec.2)}}{{
\newline > hookexp( exp(x), 9);
$$
\bigl[1,{1 \over 4},{1 \over 9},{1 \over 16},{1 \over 25},
{1 \over 36},{1 \over 49},{1 \over 64},{1 \over 81}\bigr]
$$
}} 
\smallskip\noindent
\proclaim Theorem \sec.2.
We have
$$
\sum_{T\in\setF}x^{|T|}  \prod_{h\in \setH(T)} {1\over h^2}= 
e^x. \leqno{(\sec.4)}
$$

{\it Proof}. It suffices to verify relation (\sec.2).
$$
{1\over n!}={1\over n^2} {1\over (n-1)!} + {1\over n^2} {1\over (n-2)!}. \qed
\leqno{(\sec.5)} 
$$

{\it Remark}. The number of ordered pairs of increasing labeled Fibonacci trees
on 
$\{1,2,\ldots, n\}$ 
of the same shape (i.e., the same Fibonacci tree) is equal to 
$n!$ (See, e.g., [St75, SY89]).
\medskip

\runhookexp[4pt]{18mm}{\hbox{\tt (M.\sec.3)}}{{
\newline > hookexp( exp(x+x\^{}2/2), 9);
$$
\bigl[1,{1 \over 2},{1 \over 3},{1 \over 4},{1 \over 5},{1 \over 6},
{1 \over 7},{1 \over 8},{1 \over 9}\bigr]
$$
}} 
\smallskip\noindent
\proclaim Theorem \sec.3.
We have
$$
\sum_{T\in\setF}x^{|T|}  \prod_{h\in \setH(T)} {1\over h}= 
\exp(x+x^2/2). \leqno{(\sec.6)}
$$

{\it Proof}. Let $f(x)=1+\sum_{n\geq 1}f_n x^n$ be the left-hand side of
(\sec.6). From (\sec.2) we have
$$
f_n={1\over n} f_{n-1} + {1\over n} f_{n-2}. 
\leqno{(\sec.7)} 
$$
On the other hand,
let $\sum_{n\geq 0}a_n x^n/n!$ be right-hand side of (\sec.6). We know that
$a_n$ is equal to the number of involutions of order $n$. Thus 
$$
a_n=(n-1) a_{n-2} +a_{n-1}.\leqno{(\sec.8)}
$$
Comparing (\sec.7) and (\sec.8) yields $f_n = a_n/n!$ \qed

{\it Remark}. The number of  increasing labeled Fibonacci trees on 
$[n]$
is equal to the number of involutions of order $n$ 
(see [St75, SY89]). 
\medskip

\runhookexp[4pt]{18mm}{\hbox{\tt (M.\sec.4)}}{{
\newline >  hookexp(1/(1-x), 9);
$$
\bigl[1,{1 \over 2},{1 \over 2},{1 \over 2},{1 \over 2},{1 \over 2},{1 \over 2},{1 \over 2},{1 \over 2}\bigr]
$$
}} 
\medskip\noindent
\proclaim Theorem \sec.4.
We have
$$
\sum_{T\in\setF}x^{|T|}  \prod_{h\in \setH(T), h\geq 2} {1\over 2}= 
{1\over 1-x}. \leqno{(\sec.9)}
$$

{\it Proof}. We check relation (\sec.2):
$$
1={1\over 2}  + {1\over 2} . \qed
$$

\runhookexp[4pt]{18mm}{\hbox{\tt (M.\sec.5)}}{{
\newline >  hookexp(1/(1-x)\^{}z, 6): map(factor, \%);
$$
\bigl[z,{1+z \over 4},{(z+2)(1+z) \over 9z+3},{(z+3)(z+2) \over 16z+8},
{(z+4)(z+3) \over 25z+15},{(z+5)(z+4) \over 36z+24}
\bigr]
$$
}} 
\medskip\noindent
\proclaim Theorem \sec.5.
We have
$$
\sum_{T\in\setF}x^{|T|}  \prod_{h\in \setH(T)} {(h+z-1)(h+z-2)\over h(hz+h-2)}
= 
{1\over (1-x)^z} \leqno{(\sec.10)}
$$
or
$$
\sum_{T\in\setF(n)}  \prod_{h\in \setH(T)} {(h+z-1)(h+z-2)\over h(hz+h-2)}
= 
{n+z-1\choose z-1 }. \leqno{(\sec.11)}
$$

{\it Proof}.  We check relation (\sec.2):
$$
{n+z-1\choose z-1 }= {(n+z-1)(n+z-2)\over n(nz+n-2)} (
{n+z-2\choose z-1 }+z{n+z-3\choose z-1 }
). \qed
$$

\runhookexp[4pt]{18mm}{\hbox{\tt (M.\sec.6)}}{{
\newline >  hookexp( (1-sqrt(1-4*x))/(2*x), 15);
$$
\bigl[1,1,{5 \over 3},2,{42 \over 19},{33 \over 14},{143 \over 58},
{130 \over 51},{34 \over 13},{323 \over 121},{19 \over 7},
{322 \over 117},{1150 \over 413},{45 \over 16}\bigr]
$$
\newline > guess(\%);
$$
{4(2n-3)(2n-1)\over (n+1)(5n-6)}
$$
}} 
\medskip\noindent
\proclaim Theorem \sec.6.
We have
$$
\sum_{T\in\setF}x^{|T|}  \prod_{h\in \setH(T), h\geq 2} 
{4(2h-1)(2h-3)\over (h+1)(5h-6)}
= 
{1-\sqrt{1-4x}\over 2x}
 \leqno{(\sec.12)}
$$
or
$$
\sum_{T\in\setF(n)} \prod_{h\in \setH(T), h\geq 2} 
{4(2h-1)(2h-3)\over (h+1)(5h-6)}
= 
{1\over n+1}{2n\choose n }.\leqno{(\sec.13)}
$$

{\it Proof}. We check relation (\sec.2):
$$
{1\over n+1}{2n\choose n }=  {4(2n-1)(2n-3)\over (n+1)(5n-6) }(
{1\over n}{2n-2\choose n-1 }+{1\over n-1}{2n-4\choose n-2 }
). \qed
$$

\runhookexp[4pt]{18mm}{\hbox{\tt (M.\sec.7)}}{{
\newline >  hookexp( ((1-sqrt(1-4*x))/(2*x))\^{}z, 7): 
$$
\leqalignno{
\bigl[&z,\ {3+z \over 4},\ {(z+5)(z+4) \over 9z+9},\ 
{(z+5)(z+7)(z+6) \over 8(2z+5)(z+2)},\cr
&{(z+9)(z+8)(z+7)(z+6) \over 5(5z+14)(z+5)(3+z)},\ 
{(z+10)(z+9)(z+8)(z+11) \over 36(z+4)(3+z)(z+6)}\bigr]\cr
}
$$
}} 
\medskip\noindent
\proclaim Theorem \sec.7.
We have
$$
\sum_{T\in\setF}x^{|T|}  \prod_{h\in \setH(T)} 
\rho(z;n)
= 
\Bigl({1-\sqrt{1-4x}\over 2x}\Bigr)^z,
 \leqno{(\sec.14)}
$$
where
$$\rho(z;n)={(z+2n-4)(z+2n-3)(z+2n-2)(z+2n-1)\over n(z+n-2)(z+n)(nz+4n-6)}.$$
In other words,
$$
\sum_{T\in\setF(n)} \prod_{h\in \setH(T)}  \rho(z;n)
= {z\over n!}\prod_{i=1}^{n-1}(2n-i+z). \leqno{(\sec.15)}
$$

{\it Proof}. As for the proof of Theorem \sec.6, we check relation (\sec.2).
\qed

\runhookexp[4pt]{18mm}{\hbox{\tt (M.\sec.8)}}{{
\newline >  hookexp( (1-sqrt(1-4*x\^{}2))/(2*x\^{}2)*(1+z*x), 11);
$$
\bigl[z,\; {1 \over 2z},\; {z \over 1+z^2},\; {1 \over z},\; 
{2z \over z^2+2},\; {5 \over 4z},\; {5z \over 2z^2+5},\; 
{7 \over 5z},\; {14z \over 5z^2+14},\; {3 \over 2z}\bigr]
$$
}} 
\medskip\noindent
\proclaim Theorem \sec.8.
We have
$$
\sum_{T\in\setF}x^{|T|}  \prod_{h\in \setH(T)} 
\rho(h)
= 
{1-\sqrt{1-4x^2}\over 2x^2} (1+zx),
 \leqno{(\sec.16)}
$$
where
$$
\rho(n)=\cases{
z, &if $n=1$;\cr
\noalign{\smallskip}
{2k-1\over (k+1)z}, &if $n=2k$ ($k\geq 1$);\cr
\noalign{\smallskip}
{2(2k-1)z\over (k+1)z^2+2(2k-1)}, &if $n=2k+1$ ($k\geq 1$). \cr
} \leqno{(\sec.17)}
$$
In an equivalent manner, it means that
$$
\leqalignno{
\sum_{T\in\setF(2k+1)} \prod_{h\in \setH(T)} 
\rho(h)
&= {z\over k+1}{2k\choose k };\cr
\sum_{T\in\setF(2k)} \prod_{h\in \setH(T)} 
\rho(h)
&= {1\over k+1}{2k\choose k }.\cr
}
$$

{\it Proof}. Relation (\sec.2) is being verified when $n$ is odd (resp. even)
$$
{z\over k+1}{2k\choose k }
=
{2(2k-1)z\over (k+1)z^2+2(2k-1)}
 (
{1\over k+1}{2k\choose k }
+z
{z\over k}{2k-2\choose k-1 }
)
$$
$$
(\hbox{\rm resp. }
{1\over k+1}{2k\choose k }
=
{2k-1\over (k+1)z} (
{z\over k}{2k-2\choose k-1 }
+z
{1\over k}{2k-2\choose k-1 }
)).\qed
$$

\runhookexp[4pt]{18mm}{\hbox{\tt (M.\sec.9)}}{{
\newline > hookexp( (1+x)/(1+x\^{}2), 4);
\newline {\tt Denominator is zero, no solution for n=3.}
\vskip 4pt
\newline >  hookexp( (1+x)/(1+x\^{}3), 16);
$$
[1, 0, -1, 1, 0, -1, 1, 0, -1, 1, 0, -1, 1, 0, -1, 1]
$$
}} 
\medskip\noindent
\proclaim Theorem \sec.9.
We have
$$
\sum_{T\in\setF}x^{|T|}  \prod_{h\in \setH(T)} 
\rho(h)
= 
{1+x\over 1+x^3},
 \leqno{(\sec.18)}
$$
where
$$
\rho(n)=\cases{
1, &if $n\equiv 1\mod 3$;\cr
0, &if $n\equiv 2\mod 3$;\cr
-1, &if $n\equiv 0\mod 3$.\cr
}
$$

{\it Proof}. Let
$$
\sum_{n\geq 0} f_n x^n ={1+x\over 1+x^3} = {1\over 1-x+x^2}.
$$
We have $f_{3k-1} =0$, 
$f_{3k} = (-1) f_{3k-2}$ and
$f_{3k-2} =  f_{3k-3}$. Relation (\sec.2) is then verified.\qed
\medskip

In fact, there is another generalization of Theorem \sec.4.
Consider the weight function $\rho$ that counts the leaves of Fibonacci trees. 
\runhookexp[4pt]{18mm}{\hbox{\tt (M.\sec.10)}}{{ 
\newline > [z,seq(1, i=1..6)]; 
$$
[z, 1, 1, 1, 1, 1, 1]
$$
\newline > hookgen(\%);
$$
\leqalignno{
&1+(z)x+(2z)x^2+(z^2+2z)x^3+(3z^2+2z)x^4+(z^3+5z^2+2z)x^5\cr
&+ (4z^3+7z^2+2z)x^6+(z^4+9z^3+9z^2+2z)x^7\cr
}
$$
}} 
\medskip\noindent
The above generating function corresponds to the
sequence A129710 in the on-line encyclopedia of integer sequences [Slo]
and is equal to the right-hand side of (\sec.19) below.
\proclaim Theorem \sec.10. 
We have 
$$
\sum_{\l\in\setF} x^{|\l|} \prod_{h\in \setH(\l), h =1}  z
= 
{1+(z-1)x\over 1-x-zx^2}.
\leqno{(\sec.19)}
$$


\bigskip
\medskip
{\bf Acknowledgements}.
The author wishes to thank Dominique Foata for 
helpful discussions during the preparation of this paper.

 \bigskip \bigskip


\centerline{References}

{\eightpoint

\bigskip 
\bigskip 



\livre An76|Andrews, George E.|The Theory of
Partitions|Addison-Wesley, Reading, {\oldstyle 1976}
({\sl Encyclopedia of Math. and Its Appl.,} vol.~{\bf 2})|


\article Be98|Bessenrodt, Christine|On hooks of Young diagrams|Ann. of 
Comb.|2|1998|103--110|
 
\divers BFS92| Bergeron, Fran\c cois; Flajolet,
Philippe; Salvy, Bruno|Varieties of increasing trees, {\sl Lecture 
Notes in Comput. Sci.}, {\bf 581}, Springer, Berlin, {\oldstyle 1992}|

\divers BG04|B\'eraud, Jean-Fran\c cois; Gauthier, Bruno|%
Maple package to guess closed form for a sequence of numbers,
{\tt http://www-igm.univ-mlv.fr/\~{}gauthier/GUESS.html}, see also [Kr01]|

\divers BM02|Bacher, Roland; Manivel, Laurent|Hooks and Powers of Parts in 
Partitions, 
{\sl S\'em.  Lothar. Combin.}, vol.~{\bf 47}, article B47d, {\oldstyle 2001}, 
11 pages|




\divers CY08|Chen, William Y.C.; Yang, Laura L.M.|On Postnikov's hook length
formula for binary trees, {\sl European Journal of Combinatorics}, 
in press, {\oldstyle 2008}|

\article DL08|Du, Rosena R. X.; Liu, Fu|%
$(k,m)$-Catalan Numbers and Hook Length Polynomials for Plane Trees|
European J. Combin|28|2007|1312--1321|

\article Dy72|Dyson, Freeman J.|Missed opportunities|%
Bull. Amer. Math. Soc.|78|1972|635--652|

\divers Eu83|Euler, Leonhard|The expansion of the infinite product 
$(1-x)(1-xx)(1-x^3)(1-x^4)(1-x^5)(1-x^6)$ etc. into a single series, 
{\sl English translation from the Latin by Jordan Bell} 
on {\it arXiv:math.HO/0411454}|



\article FH01|Foata, Dominique; Han, Guo-Niu|Arbres minimax et polyn\^omes 
d'Andr\'e|Advances in Appl. Math.|27|2001|367--389|

\divers FS73|Foata, Dominique; Sch\"utzenberger, Marcel-Paul|Nombres 
d'Euler et permutations alternantes, {\sl A survey of Combinatorial Theory} 
J.N. Srivastava et al., eds., p. 173--187. Amsterdam, North-Holland, 
{\oldstyle 1973}| 

\article FS74|Foata, Dominique; Strehl, Volker|Rearrangements of 
the symmetric group and enumerative properties of the tangent and secant 
numbers|Math. Zeitschrift|137|1974|257--264| 

\article FRT54|Frame, J. Sutherland; Robinson, Gilbert de Beauregard;        
Thrall, Robert M.|The hook graphs of the symmetric groups|Canadian 
J. Math.|6|1954|316--324|

\divers Ga01|Garvan, Frank|%
A $q$-product Tutorial for a $q$-series Maple Package,
{\sl The Andrews Festschrift. Seventeen Papers on Classical Number Theory 
and Combinatorics}, D. Foata, G.-N. Han eds., Springer-Verlag, 
Berlin Heidelberg, {\oldstyle 2001}, pp. 111-138.
{\sl Sem. Lothar. Combin.} Art. B42d, 27 pp|

\article GKS90|Garvan, Frank; Kim, Dongsu; Stanton, Dennis|Cranks and 
$t$-cores|Invent. Math.|101|1990|1--17|

\article GNW79|Greene, Curtis; Nijenhuis, Albert;
Wilf, Herbert S.|A probabilistic proof of a formula for the number of     
Young tableaux of a given shape|Adv. in Math.|31|1979|104--109|

\divers GS06|Gessel, Ira M.; Seo, Seunghyun|%
A refinement of Cayley's formula for trees,
{\it arXiv:math.CO/0507497}|

\article GV85|Gessel, Ira; Viennot, Gerard|%
Binomial determinants, paths, and hook length formulae|Adv. in  
Math.|58|1985|300--321|
  
\divers Ha08a|Han, Guo-Niu|An explicit expansion formula for the powers of 
the Euler Product in terms of partition hook lengths, 
{\sl arXiv:0804.1849v2, Math.CO}, 35 pages, {\oldstyle 2008}|

\divers Ha08b|Han, Guo-Niu|New hook length formulas for binary trees,
{\sl arXiv:0804.3638, Math.CO}, 4 pages, {\oldstyle 2008}|

\divers Ha08c|Han, Guo-Niu|Yet another generalization of 
Postnikov's hook length formula for binary trees,
{\sl arXiv:0804.4268v1, math.CO} , 4 pages, {\oldstyle 2008}|

\divers Ha08d|Han, Guo-Niu|\quad Some conjectures and open problems
about partition hook length, {\sl in preparation}, 12 pages, {\oldstyle 2008}|

\divers Ha08e|Han, Guo-Niu|The Nekrasov-Okounkov hook length formula: 
refinement, elementary proof, extension and applications,
{\sl in preparation}, 28 pages, {\oldstyle 2008}|


\livre JK81|James, Gordon; Kerber, Adalbert|The
representation theory of the symmetric group|Encyclopedia
of Mathematics and its Applications, {\bf 16}. Addison-Wesley Publishing, 
Reading, MA, {\oldstyle 1981}|

\article JS89|Joichi, James T.; Stanton, Dennis|An
involution for Jacobi's identity|Discrete Math.|73|1989|261--271|


\article Kn70|Knuth, Donald E.|Permutations, matrices, and generalized Young
tableaux|Pacific J. Math.|34|1970|709-727|

\livre Kn98a|Knuth, Donald E.|The Art of Computer Programming,  {\bf vol.~1}, 
Fundamental Algorithms, 3rd ed.|Addison Wesley Longman,  {\oldstyle 1997}|

\livre Kn98b|Knuth, Donald E.|The Art of Computer Programming,  {\bf vol.~3}, 
Sorting and Searching, 2nd ed.|Addison Wesley Longman,  {\oldstyle 1998}|



\divers Kr01|Krattenthaler, Christian|%
\quad RATE - A Mathematica guessing machine, 
{\tt http:// igd.univ-lyon1.fr/\~{}kratt/rate/rate.html}, see also [BG04]|

\article Kr99|Krattenthaler, Christian|\ Another involution
principle-free bijective proof of  
Stanley's hook-content formula|J.      
Combin. Theory Ser. A|88|1999|66--92| 


\livre La01|Lascoux, Alain|Symmetric Functions and Combinatorial Operators on 
Polynomials|CBMS Regional Conference Series in Mathematics, Number 99, 
{\oldstyle 2001}|


\article Ma72|Macdonald, Ian G.|\ Affine root systems and 
Dedekind's $\eta $-function|Invent. Math.|15|1972|91--143|

\livre Ma95|Macdonald, Ian G.|Symmetric Functions and Hall Polynomials|
Second Edition, Clarendon Press, Oxford, {\oldstyle 1995}|



\article MY07|Moon, J. W.; Yang, Laura L. M.|%
Postnikov identities and Seo's formulas|Bull. Inst. Combin. Appl.|%
49|2007|21--31|

\divers NO06|Nekrasov, Nikita A.; Okounkov, Andrei|Seiberg-Witten theory and 
random partitions. The unity of  mathematics, 525--596, {\sl Progr. Math.}, 
{\bf 244}, Birkhaeuser Boston. {\oldstyle 2006}. (See also 
{\sl arXiv:hep-th/0306238v2}, 90 pages, {\oldstyle 2003})|

\article NPS97|Novelli, Jean-Christophe; Pak, Igor;
Stoyanovskii, Alexander V.|A direct bijective proof of the hook-length
formula|Discrete Math. Theor. Comput. Sci.|1|1997|53--67|

\divers Po04|Postnikov, Alexander|Permutohedra, associahedra, and beyond,
{\it arXiv:math. CO/0507163}, {\oldstyle 2004}|

\divers Pr96|Prodinger, Helmut|
A Note on the Distribution of the Three Types of Nodes in Uniform Binary Trees,
with comments by Christian Krattenthaler, 
Guo-Niu Han and G\"unter Rote, 
{\it S\'eminaire Lotharingien de Combinatoire}, Article B38b,  
{\oldstyle 1996}, 5 pp|


\article RW83|Remmel, Jeffrey B.; Whitney, Roger A|
bijective proof of the hook formula for the number of column strict tableaux
with bounded entries|European J. Combin.|4|1983|45--63|

\divers Sch76|Sch\"utzenberger, Marcel-Paul|La corres\-pondance de Robinson,
dans ``Combinatoire et Repr\'esentation du 
Groupe Sym\'etrique", {\sl Lecture Notes in Mathematics}, Springer-Verlag, 
vol. {\bf 579}, {\oldstyle 1976}, p. 59--113|

\divers Se08|Seo, Seunghyun|%
A combinatorial proof of Postnikov's identity and a generalized enumeration 
of labeled trees,
{\it arXiv:math.CO/0409323}|


\divers Slo|Sloane, Neil; {\it al.}|
The On-Line Encyclopedia of Integer Sequences, {
\tt http:// www.research.att.com/\char126njas/sequences/}|

\article St75|Stanley, Richard P.|The Fibonacci lattice|%
Fibonacci Quart.|13|1975|215--232|

\article St76|Stanley, Richard P.|Theory and application of 
plane partitions (II)|Studies in Appl. Math.|50|1971|259--279|

\livre St97|Stanley, Richard P.|Enumerative Combinatorics, vol. 1|
Cambridge university press, {\oldstyle 1997}|

\livre St99|Stanley, Richard P.|Enumerative Combinatorics, vol. 2|
Cambridge university press, {\oldstyle 1999}|



\article SY89|Sagan, Bruce E.; Yeh, Yeong Nan|Probabilistic algorithms 
for trees|Fibonacci Quart.|27|1989|201--208|

\divers Vi81|Viennot, G\'erard|\quad Interpr\'etations combinatoires des 
nombres d'Euler et de Genocchi, {\it S\'eminaire de Th\'eorie des Nombres},
Bordeaux, 94 pages, {\oldstyle 1981}|





\article Ze84|Zeilberger, Doron|A short hook-lengths bijection
inspired by the Greene-Nijen\-huis-Wilf proof|Discrete 
Math.|51|1984|101--108|

\bigskip

\irmaaddress
}
\vfill\eject

\end